\documentclass{article}

\usepackage{PRIMEarxiv}

\usepackage[utf8]{inputenc} % allow utf-8 input
\usepackage[T1]{fontenc}    % use 8-bit T1 fonts
\usepackage{hyperref}       % hyperlinks
\usepackage{url}            % simple URL typesetting
\usepackage{booktabs}       % professional-quality tables
\usepackage{amsfonts}       % blackboard math symbols
\usepackage{nicefrac}       % compact symbols for 1/2, etc.
\usepackage{microtype}      % microtypography
\usepackage{lipsum}
\usepackage{fancyhdr}       % header
\usepackage{graphicx}       % graphics
\graphicspath{{media/}}     % organize your images and other figures under media/ folder

\usepackage{caption}
\usepackage{amsfonts,amssymb,amsthm}       % 

\usepackage{hyperref}
\usepackage{url}
\usepackage{graphicx}
\usepackage{booktabs}
\usepackage{natbib}
\usepackage{mathtools}
\usepackage{algorithm}
\usepackage{algpseudocode}

\newtheorem{theorem}{Theorem}

\newtheorem{lemm}{Lemma}

\newcommand{\bfm}[1]{\boldsymbol{#1}}
  
\title{Partially-Supervised Neural Network Model For Multiparametric Quadratic Programming
}

\author{
  Fuat Can Beylunioğlu, Mehrdad Pirnia, P. Robert Duimering \\
  Management Science and Engineering \\
  University of Waterloo \\
  Waterloo\\
  \texttt{\{fcbeylun, mpirnia,rduimering\}@uwaterloo.ca} \\
}

\begin{document}
\maketitle

\begin{abstract}
Neural Networks (NN) with ReLU activation functions are used to model multiparametric quadratic optimization problems (mp-QP) in diverse engineering applications. Researchers have suggested leveraging the piecewise affine property of deep NN models to solve mp-QP with linear constraints, which also exhibit piecewise affine behaviour. However, traditional deep NN applications to mp-QP fall short of providing optimal and feasible predictions, even when trained on large datasets. This study proposes a partially-supervised NN (PSNN) architecture that directly represents the mathematical structure of the global
solution function. In contrast to generic NN training approaches, the proposed PSNN method derives a large proportion of model weights directly from the mathematical properties of the optimization problem, producing more accurate solutions despite significantly smaller training data sets.  Many energy management problems are formulated as QP, so we apply the proposed approach to energy systems (specifically DC optimal power flow) to demonstrate proof of concept. Model performance in terms of solution accuracy and speed of predictions was compared against a commercial solver and a generic Deep NN model based on classical training. Results show KKT sufficient conditions for PSNN consistently outperform generic NN architectures with classical training using far less data, including when tested on extreme, out-of-training distribution test data. Given its speed advantages over traditional solvers, the PSNN model can quickly produce optimal and feasible solutions within a second for millions of input parameters sampled from a distribution of stochastic demands and renewable generator dispatches, which can be used for simulations and long term planning.
\end{abstract}

% keywords can be removed

\section{Introduction}

There has been increasing interest in using Neural Networks (NN) to predict the solutions to complex nonlinear optimization problems in energy management, chemistry, control theory, and other domains \citep{nellikkath2022physics,karg2020efficient,katz2020integrating}. Effectively, this goal amounts to estimating a solution function that predicts optimal solutions based on a given set of input parameters, such as the right-hand-side (RHS) vector and cost coefficients of an optimization model. Most applications treat the NN model as a black-box and employ standard training methods using datasets of input-output pairs, representing particular values of system parameters and corresponding optimal solutions. Most also ignore the mathematical structure of the underlying function and utilize a generic deep NN (DNN) architecture to approximate solutions. Such modeling approaches necessitate large, computationally expensive training datasets to achieve satisfactory performance, yet cannot guarantee the feasibility or optimality of results.

% Through multi-parametric programming, one can obtain the optimization variables of the
% problem as a function of the bounded uncertain parameters, and
% the regions (in the space of the parameters) where these functions are valid.

% The goal is to train a model that can predict optimal solutions to optimization problems without going through expensive iterative calculations of traditional solvers.

% the function producing optimal solutions for any given input parameters with a multiparameter programming perspective 
% The mathematical properties of the optimal solutions with respect to the changes in the model parameters is well studied in the multiparametric programming literature. 
% Many optimization problems involve parameters that are
% unknown, either because they cannot be measured, or because
% they represent information about the future (e.g. future state of
% a system, future disturbance, future demand). Multi-parametric
% programming is a technique for the solution of such class of
% uncertain optimization problems.

% Through multi-parametric programming, one can obtain the optimization variables of the problem as a function of the bounded uncertain parameters, and the regions (in the space of the parameters) where these functions are valid.

% Multiparametric programming is an approach for obtaining the optimal solutions to the problem as a function of feasible problem parameters, and identify the regions of parameter space where these functions are valid.
The multiparametric programming literature similarly seeks to represent the solution function to optimization problems, using algebraic methods to characterize optimal solutions as a function of feasible problem parameters\footnote{The function mapping parameters to optimal primal solutions is sometimes referred to as an optimizer in the multiparametric programming literature. We use the term solution function to refer to the mapping between parameters to both primal and dual solutions.} and to identify the regions of parameter space where these functions are valid \citep{pistikopoulos2020multi}. Recent work has attempted to integrate multiparametric programming with DNNs, based on recognition of their similar underlying mathematical characteristics \citep{karg2020efficient,katz2020integrating}. For multiparametric Linear Programs (mp-LP) and Quadratic Programs (mp-QP) with linear constraints, the optimal primal and dual solutions are piecewise linear (PWL) functions of the parameters \citep{pistikopoulos2020multi}. Specifically, each set of constraints binding at a solution corresponds with a subregion of the feasible domain, called a \textit{critical region}, where solutions change linearly with respect to changes in the RHS parameters, and solutions change piecewise linearly between adjacent critical regions. Similarly, a NN with ReLU activation functions is a general form of PWL function with trainable weights and biases \citep{montufar2014number}. In theory, therefore, for any mp-LP or mp-QP one should be able to find an optimal set of weights and biases such that a NN can exactly represent the corresponding solution function \citep{karg2020efficient}. 

In practice, however, DNN applications of multiparametric programming face significant training challenges to guarantee global solution optimality and feasibility. First, obtaining the optimal NN weights and biases is difficult due to nonconvexity of the training loss function, whereby multiple sets of weights and biases can result in different local minima for a given dataset. Second, NN prediction accuracy is highly dependent on the training set, which makes it hard to generalize predictions outside the training distribution. Overcoming this problem to find the global solution function requires sampling very large datasets that cover all critical regions of the parameters, which is impractical for problems with high dimensions. Specifically, for an $n$ dimensional domain, constructing a dataset that includes only boundary points for each feasible region requires sampling $2^n$ points, which would still be insufficient because additional interior points are needed within each region for good estimation.

To address these concerns, this paper first considers some of the challenges involved in training a NN model to represent a PWL function without error and then proposes a partially-supervised NN (PSNN) architecture and training procedure for mp-QP problems that directly reflects the mathematical structure of the global solution function. In our analysis, we show that for a NN model to represent the target function exactly, two conditions must hold: (i) the model layer size must be equal to the number of linear segments in the target PWL function, and (ii) a minimal number of training data points must be sampled on each linear segment of the target function. Our proposed model and training approach satisfy these conditions, resulting in solution predictions that satisfy KKT sufficient optimality and feasibility requirements while dramatically reducing dataset demands compared to traditional DNN training methods. In general, while the solution function to mp-QP for uncertain parameters is piecewise linear, the slope changes in the solution function can be decomposed into two components: piecewise linear changes in the value of dual variables associated with inequality constraints, and the remaining terms that form a linear function. The proposed PSNN implements these components separately in a single NN architecture, whose calculations can be easily parallelized using GPU to predict solutions for a large set of input parameters. The proposed PSNN approach was used to model DC optimal power flow (DC-OPF) problems in the domain of electricity power management for proof-of-concept support.

Unlike black-box NN models, the proposed PSNN for QP is constructed partially supervised by deriving a large part of the NN weights directly from the problem. These coefficients are calculated prior to the training by expanding the Lagrangian function with each possible combination of binding constraints, calculating partial derivatives and inverting the resulting Jacobian matrices of the derivatives of the Lagrangian function. Theoretically, the prior computation of coefficients can be infeasible as the number of constraint combinations increases exponentially. However, MP literature discusses that under certain conditions, when moved from one critical region to an adjacent one, either one constraint is added or dropped from the set of binding constraints. Therefore, potential sets of binding constraints linearly increase, which can be discovered by solving the problem for a set of input system parameters sampled equidistantly towards a certain direction. This procedure is used only to filter out unused cases and is run only about 500 times for even the largest empirical test systems investigated. Once the set of possible binding constraints is discovered, this information is used to construct datasets with minimal computational costs via linear operations and without using solvers in a loop. 

%Unlike black-box NN models, the proposed PSNN for QP is largely constructed unsupervised by deriving NN weights directly from the problem. The coefficients are calculated prior to the training by expanding the Lagrangian function with each possible combination of binding constraints, calculating partial derivatives and inverting the resulting Jacobian matrices of the derivatives of the Lagrangian function. Theoretically, the prior computation of coefficients can be infeasible for large problems, as the number of constraint combinations increases exponentially. Nevertheless, in practical applications, many constraint combinations are not possible due to characteristics of the problem, such as the topology of the electricity grid in energy management. To scale down the set of all possible binding constraint combinations, we use a commercial solver to solve the problem for a set of input system parameters sampled equidistantly from minimum to maximum feasible input demand for a number of arbitrarily chosen parameters. Note that this procedure is used only to filter out unused cases and is run only about 500 times for even the largest empirical test systems investigated. Once the set of possible binding constraints is discovered, this information is used to construct datasets with minimal computational costs without using solvers in a loop.

Our paper differs from prior literature in the following ways. Firstly, we focus on exact representation of the solution function and addresses the challenges. We propose a training procedure that share information with the true solution function and learns to align with this function. In contrast, existing studies typically train black-box approaches that terminate training once validation loss no longer decreases, without providing insight into how to bridge the gap between theoretical potential and actual performance. This oversight is particularly addresses a gap in mp-QP literature, which posits that NN is a general form of piecewise-linear solution functions but falls short of providing a clear roadmap for achieving this goal. Moreover, the emphasis on representing the solution function results in models that generalize well beyond the training distribution, outperforming existing approaches whose performances are limited to the training set range.

% Our paper differs from prior literature in the following ways. Our study emphasizes representing the solution function exactly, addresses issues on this goal and proposes an approach to achieve this goal. This contrasts with the existing studies that train that terminate after validation loss no longer reduces. This addresses one fundamental gap mp-QP literature  that proposes NNs can naturally represent the PWL solution function but does not provide further information to how to achieve this goal. Thanks to this alignment, our trained models can generalize outside the training distribution well unlike the existing approaches that further process the calculated results to reduce prediction errors.

% Therefore, there is no clear state-of-the-art for us to compare our method. However, most studies report larger errors than our results reported in this paper, which are measured with  optimality gaps of 1E-6 to 1E-10 and squared KKT violations less than 1E-8. The MP literature focuses primarily on methodology and performances are not explicitly reported \cite{karg2020efficient, katz2020integrating}. On the other hand, PINN model trained to solve DC-OPF problem reports 0.12\% to 3.11\% KKT violations \cite{nellikkath2021DCPINN}, while the NeuralDeco model reports 0.5\% optimality gap in their solutions \cite{chen2022learning}. DC3 model trained to solve quadratic problems reports 6.219 (0.098\% in percentage) optimality gap. Note that this model further processes the results using gradient-descent operations when the predictions do not satisfy the feasibility conditions. 

This study contributes to the literature as follows:
\begin{itemize}
	\item We show why classical DNN training methods cannot yield a NN model that represents the solution function to mp-QP problem with linear constraints without error.
	% 		\item A scalable dataset construction method is proposed that bypasses running solvers in a loop, which generate data points with random noise.
	
	\item We propose a partially supervised NN model that can provide optimal solutions to QP with linear constraints with uncertain parameters on the RHS of equality constraints. The model is based on an explainable NN architecture that directly aligns with the piecewise linearity of the underlying optimization problem, and is supported by a training approach that decomposes the learning process to subproblems to overcome the training challenges.

	\item The PSNN modeling approach is applied to DC-OPF problems with generator capacity limits for proof-of-concept support. The resulting models provide precise predictions that satisfy KKT optimality and feasibility conditions, and are used rapidly to create a distribution of optimal and feasible solutions to a large dataset sampled from the empirical distribution of uncertain demand and renewable generator dispatches.
	
\end{itemize}

The paper is organized as follows: Section \ref{sec:lit} discusses related studies in the literature. Section \ref{sec:background} presents background on NNs as universal function approximators and the piecewise linearity of NN with ReLU activation, and considers trade-offs between model complexity and dataset requirements to represent the target function exactly. Section \ref{sec:method} details our methodology, discusses how QP solutions form a PWL function that can be decomposed to linear and piecewise parts, and outlines our partially supervised NN architecture and training procedure. Section \ref{sec:dcsolver} presents numerical results, and sections \ref{sec:summary} and \ref{sec:ssnn_limitations} discuss conclusions, limitations and future work.

\section{Related Work}\label{sec:lit}

% Studies first proposed using NN for solving system of equations in the 1990s, when training NNs with high computation power was not available. Early studies \cite{lee1990neural,lagaris2000neural} leveraged the approximation power of NNs to solve ordinary or partial differential equations (PDEs). Later researchers applied NN models to atomic and molecular physics \cite{caetano2011using} and other science applications  \cite{parisi2003solving,piscopo2019solving}. Despite their advantages as nonlinear solvers, training a NN to solve a system of equations is challenging because finding the best model parameters from a large parameter space is data-greedy and computationally challenging for optimization algorithms. 

\subsection{Surrogate Models}
Using NNs for surrogate modeling of optimization problems is an established field with various applications in engineering such as control \cite{karg2020efficient}, chemistry \cite{katz2020integrating}, and power systems management \citep{nellikkath2022physics,fioretto2020predicting,lotfi2022constraint}. The approach leverages DNNs as universal function approximators to estimate the nonlinear relationship between problem parameters and optimal solutions. Despite the potential advantages of DNNs as nonlinear solvers, guaranteeing feasibility and optimality is challenging, and finding the best model parameters from a large parameter space is data-greedy with significant computational training burdens.  

% However, guaranteeing optimality and feasibility of solutions for all inputs poses a challenge due to the nonlinearity of the solution function.

To improve feasibility, \cite{zamzam2020learning} and \cite{chen2023end} proposed two-stage models that first predict optimal solutions with a NN and then post-process the solutions to enforce feasibility.  \cite{lotfi2022constraint} proposed a NN architecture that enforces optimal dispatch generator limits using a sigmoid activation function in the output layer. \cite{fioretto2020predicting} processed optimal solutions through sequentially connected sub-networks reflecting problem characteristics. Physics-informed Neural Networks (PINNs) use DNNs that are also trained to satisfy KKT optimality conditions, resulting in reduced dataset requirements, more accurate predictions with less violations, and more generalizability power \citep{nellikkath2022physics}.  

Other studies in the OPF literature are concerned with using ML to improve execution time of solvers. In a series of studies (e.g. \cite{misra2021learning,deka2019learning}), researchers observed that it is possible to reduce the dimensions of the problem by predicting the active constraints. Specifically, the solution must always satisfy the equality constraints and predicting the binding inequality constraints reduces the size of the problem space. For this purpose, different types of classifiers were trained to find the solution faster \cite{misra2021learning,deka2019learning,ng2018statistical}. Others \cite{baker2018joint} proposed faster solutions to the AC-OPF problem by removing inactive constraints to modify the feasibility space and downsize the problem.

% Specifically, the solution must always satisfy the equality constraints and predicting the binding inequality constraints reduces the size of the problem space. 

% The ability to use neural network models as surrogates in an optimization formulation is contingent on utilizing a strategy to obtain a quality optimal solution. Current literature demonstrates an interest in obtaining ‘good’ solutions to optimization formulations incorporating deep learning models. Schweidtmann et al. [10] proposed a global optimization strategy based on McCormick relaxations to identify the global minimum to optimization problems involving ANNs. Pfrommer et al. [15] utilized a stochastic genetic algorithm to find the minimum for a process involving textile draping where a neural network was utilized as a surrogate model. Nagat et al. [16] developed a surrogate neural network model for a fermentation process, and optimal operating conditions were identified using a genetic algorithm.

\subsection{Multiparametric Programming}
Multiparametric programming became more popular in the 2000s when researchers showed that model predictive control laws can be expressed as multiparametric programming solutions \cite{bemporad2000explicit}, specifically in the form of a piecewise linear function of the parameter set. Studies proposed a geometric interpretation that considers the critical regions of mp-LP solutions as polyhedra, where the parameter space is explored by moving the parameter set among critical regions \cite{borrelli2003geometric}.  However, it was later shown that this approach does not guarantee discovery of all critical regions \cite{tondel2003algorithm}.  The geometric approach was later generalized to solve mp-QP by Bemporad et al. \cite{bemporad2002explicit} as it shares similar characteristics with mp-LP. In contrast to mp-LP, which relies on the simplex algorithm, multiparametric QP was theorized using the Basic Sensitivity Theorem proposed by Fiacco \cite{fiacco1976sensitivity}. It was shown that solutions to mp-QP are PWL functions and the optimal objective function is a piecewise quadratic function of the parameter set \cite{bemporad2000explicit, pistikopoulos2020multi}. 

The geometric approach involves finding critical regions by solving the problem for an initial parameter set, $\bfm\theta_0$, identifying active sets and obtaining the parametric solution corresponding to the critical region where the solution is valid. Subsequently, the remaining critical regions are explored by moving outside the initial critical region. The original study by \cite{bemporad2000explicit} proposes constraint reversal by successively reversing signs of each constraint, i.e., from $A_ix\leq b_i$ to $A_ix> b_i$. However, this approach creates too many artificial cuts leading to scalability challenges for larger problems. Two main modifications were proposed to overcome scalability. First is an active set inference approach \cite{tondel2003algorithm} that adds or drops one constraint at a time from the active constraint set to move to an adjacent region. The second is a variable step size approach \cite{baotic2003new} that chooses the facets defining the critical region, finding its centre and moving outside to discover an adjacent critical region. However, both methods fail to explore the entire feasible parameter set when facet-to-facet property does not hold.

The critical regions of an mp-QP problem are each defined by a unique corresponding set of active constraints, so the feasible space can be exhaustively explored by enumerating all combinations of constraints and removing infeasible combinations. However, as the number of constraint combinations grows exponentially, this approach does not scale up to higher dimensional problems. The number of potential combinations can be reduced when an infeasible set is discovered \cite{gupta2011novel} since a power of an infeasible set is also infeasible, so branch and bound type algorithms have been proposed to discover the critical regions \cite{feller2012combinatorial,feller2013explicit}.

The mixed integer versions of LP and QP problems can be viewed as combinations of multiple mp-LP and mp-QP problems. Therefore, strategies proposed to solve mp-MILP and mp-MIQP problems include enumerating all combinations of integer variables, solving the corresponding LP or QP problems, and combining solutions by comparison. To make the solution tractable, researchers again proposed branch and bound  \cite{acevedo1999algorithm} and decomposition strategies  \cite{jia2006uncertainty}. 

The properties of critical regions are different for mp-MIQP and mp-MILP. Whereas the critical regions for mp-MILP are polytopes, obtaining optimal solutions by comparing multiple mp-QP solutions can result in quadratically constrained critical regions. Maintaining critical regions in the form of polyhedra offers computational advantages, so researchers have applied McCormick relaxations for linearizing the objective function or critical regions \cite{oberdieck2015explicit}. Other applications of multiparametric programming to optimization, including multiparametric nonlinear programming, are outside the scope of this study (see \cite{oberdieck2016multi} for further review.)

Examining mathematical properties of mp-LP and mp-QP problems under a multiparametric setting led researchers to recognize their similarity to the PWL characteristics of DNNs.  \cite{katz2020integrating} attempted to integrate PWL behaviour of the problem solution analytically. \cite{karg2020efficient} examined the functional similarities between DNNs and model predictive control (MPC) laws, and trained a model to estimate MPC solutions for a time variant system. Similarly, \cite{huo2022integrating} proposed an approach to integrate MPC laws of a MILP problem for a microgrid system. However, these approaches employ black-box models that do not exactly represent the solution function. Moreover, none of these works compare optimality and feasibility between their models and traditional solver benchmarks or DNN based approaches. This paper proposes a custom NN architecture that closely aligns with the PWL structure of the global solution function to improve accuracy and generalizability of mp-QP predictions.

%The outline of this paper is as follows: The next section illustrates the challenge of training a black-box NN model to represent the solution function exactly. Section \ref{sec:method} discusses how QP solutions can be decomposed to linear and PWL parts and presents proofs to support the argument. Later, the section details our proposed NN architecture that provides exact solution to the linear and PWL functions and combines both in one NN flow. Section \ref{sec:dcsolver}  details model training and compares our model performance against alternative approaches on test cases. Section \ref{sec:summary} discusses some of the challenges that prevent the model from scaling up to larger problems and concludes the paper.

\section{Background}\label{sec:background}

\subsection{Neural Networks}\label{sec:NNs}
Neural Networks are universal function approximators, defined as composite functions that can be expressed as an indirect mapping between inputs $\textbf{x}$ and outputs $\textbf{y}$, $f: \textbf{x} \rightarrow  \textbf{y}$. The basic type of NN, known as shallow feed forward neural network, consists of an input and output layer, connected by a hidden layer, $ h $. Loosely influenced by neurons and synapses of the brain, the value of each neuron is determined by the weighted sum of the previous neurons' values, which is called \textit{input current} denoted by $ z_i $ for node $i$. However, the neuron's value is not always defined as its input current but usually transformed by an activation function. For example, the rectified linear unit (ReLU) activation function is a threshold function that is activated only if the input current value is positive, expressed mathematically as follows:
\begin{equation}
	z_i^1 = \sum_{j=1}^{d_x}W_{ji}^0x_j + b^0_i, \quad h_i = \sigma(z_i^1), 
\end{equation}
where $z_i^1$ is input current of node $i$ in the hidden layer, $ W^0_{ji} $ is the weight of the edge connecting node $ j $ of one layer to node $ i $ of the next layer (here input layer to hidden layer), $ b_i^0 $ is the bias term adjusting the calculated value for neuron $i$, and $ \sigma(z_i^1)=ReLU(z_i^1) $ is an element-wise activation function defined as
\begin{equation}
	ReLU(z_i) = \begin{cases}
		z_i& \text{ if } z_i\geq 0,\\
		0 & \text{ otherwise. }
	\end{cases}
\end{equation}
The output layer is defined by the previous layer, $ h(\textbf{x}) $, in a similar fashion as
\begin{equation}
	z^2_k = \sum_{i=1}^{d_h} W^1_{ik}h_i + b_k^1, \quad y_k = \sigma(z^2_k). 
\end{equation}
Here, the activation function $ \sigma(\cdot) $ is not necessarily the same as the one activating the previous layer. For simplicity, a generic notation $\sigma(\cdot)$ will be used to denote all activation functions in this paper.

\begin{figure}
	\centering
	\begin{minipage}{.4\textwidth}
		\centering
		\vspace{5pt}
		\includegraphics[width=.8\linewidth]{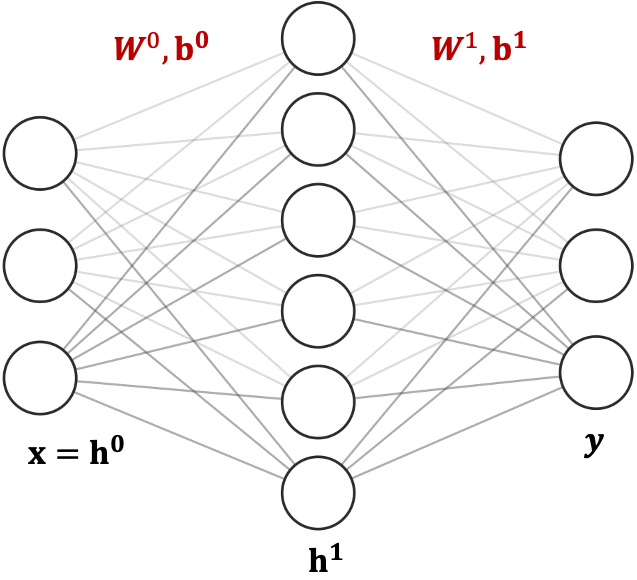}
	\end{minipage}%
	\begin{minipage}{.55\textwidth}
		\centering
		\vspace{2pt}
		\includegraphics[width=.8\linewidth]{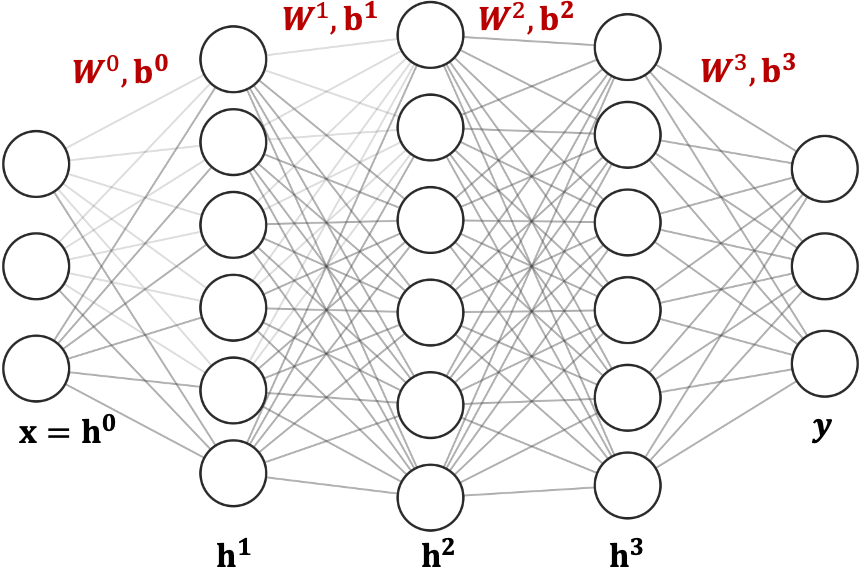}
		
	\end{minipage}
	\caption{Illustration of Shallow (left) and Deep (right) NN models}
	\label{fig:NN_illust}
\end{figure}

The above example expresses a \textit{shallow network} with one hidden layer shown as the diagram on the left in Figure \ref{fig:NN_illust}. In the case of a \textit{deep neural network} with two or more hidden layers (the diagram on the right in Figure \ref{fig:NN_illust}b), the output value can be calculated iteratively as follows:
\begin{align} \label{eq:deepNN}
	\textbf{h}^0 &= \textbf{x}^T,\notag \\
	\textbf{h}^{i} &= \sigma(\textbf{W}^{i-1T} \textbf{h}^{i-1T}+\textbf{b}^{i-1}), \quad  \text{ for } i \in [1,\dots,n_h],\notag \\
	\textbf{y}  &= \sigma(\textbf{h}^{n_hT}),
\end{align}
where $ n_h $ is the number of hidden layers, the superscript $T$ denotes the transpose operation and $\textbf{W}^{iT} = (\textbf{W}^{i})^T$ is used to simplify the notation, $ \textbf{x} = [x_1,\dots,x_{d_x}] $, $ \textbf{h}^i = [h_1,\dots,h_{d_{h^i}}] $, $ \textbf{y} = [y_1,\dots,y_{d_y}] $, $\textbf{x}\in \mathbb{R}^{d_d\times d_x}$, $\textbf{h}^i\in \mathbb{R}^{d_d\times d_{h^i}}$, $\textbf{y}\in \mathbb{R}^{d_d\times d_y}$. Here, $ \textbf{W}^i\in \mathbb{R}^{n\times m}$ and $\textbf{b}^i\in \mathbb{R}^{n\times 1}$ are parameters known as weight and bias terms that are optimized during NN training, and generically defined as,
\begin{align}
	\textbf{W}^i = 
	\begin{bmatrix}
		W_{11}^i & \dots & W_{1m}^i\\
		\vdots & \ddots & \vdots \\
		W_{n1}^i & \dots & W_{nm}^i\\
	\end{bmatrix}, \quad 
	\textbf{b}^i = 
	\begin{bmatrix}
		b_1^i \\ \vdots \\b_n^i
	\end{bmatrix}.
\end{align}
For a shallow model where $n_h=1$, $ \textbf{W}^1 \in \mathbb{R}^{d_x\times d_h} $, $ \textbf{b}^1 \in \mathbb{R}^{d_h} $,$ \textbf{W}^2 \in \mathbb{R}^{d_h\times d_y} $ and $ \textbf{b}^2 \in \mathbb{R}^{d_y} $. 

%In this study, NNs are often compared with basic linear models, $ L(\textbf{x};\theta) $, or linear models activated by some activation function. These models can be seen as equivalent to a network with no hidden layer. For simplicity, these models will often be referred to as NN with $n_h=0$ in the rest of the study. The linear model can be mathematically expressed as
%\begin{equation}
%	L(\textbf{x};\theta) = \sigma(\textbf{W}\textbf{x} + \textbf{b}).
%\end{equation}
%Because the model is linear, the activation is an identity function, i.e. $ \sigma(\textbf{x}) = \textbf{x} $ but could be activated by different functions to incorporate nonlinearity as well. The linear model is illustrated in Figure \ref{fig:fflinear}. 
%
%Also, the performance of various types of networks will be compared throughout the study. For the sake of clarity, a neural network will be denoted as,
%$f^\mathcal{A}(\textbf{x}; \theta),$
%where $ \theta $ is the parameter vector, $ \mathcal{A} $ is the architecture that contains the information such as $ d_h, n_h, \sigma(x) $. For example, $ L(\textbf{x}; \theta) $ can also be denoted as $ f^\mathcal{A}(\textbf{x};\theta) $ where $ \mathcal{A} = \{d_h=0,n_h=0, \sigma(x) = x\} $.
%
%Finally, it is worth mentioning that there are various other architectures of neural network besides Feed Forward NN. Our study only focuses on Feed Forward NN types and its variations and the other types such as RNN, CNN and networks with attention mechanisms are out of the scope of this work.

\subsection{The Tradeoff Between Dataset Size and Model Complexity}\label{sec:armsrace}
As our goal is to train a NN model that represents the solution function with no error, we conducted preliminary analysis to determine whether a NN can be trained to represent a target function via standard training. Specifically, we investigated whether it is possible to train a DNN on data sampled from a target function. 

In this section, we demonstrate why traditional NN training practices are not effective for finding optimal model weights and biases that can result in an exact representation of the target function, i.e., producing zero train and test error. Specifically, our analysis shows that there is a minimal model complexity and dataset size required to accurately estimate the solution function exactly. However, if a more complex model than minimal complexity is trained, it becomes practically impossible to represent the target function exactly and the model error can only be reduced with larger dataset size.

\begin{figure*}
	\centering
	\begin{minipage}{.45\textwidth}
		\centering
		\includegraphics[width=.95\linewidth]{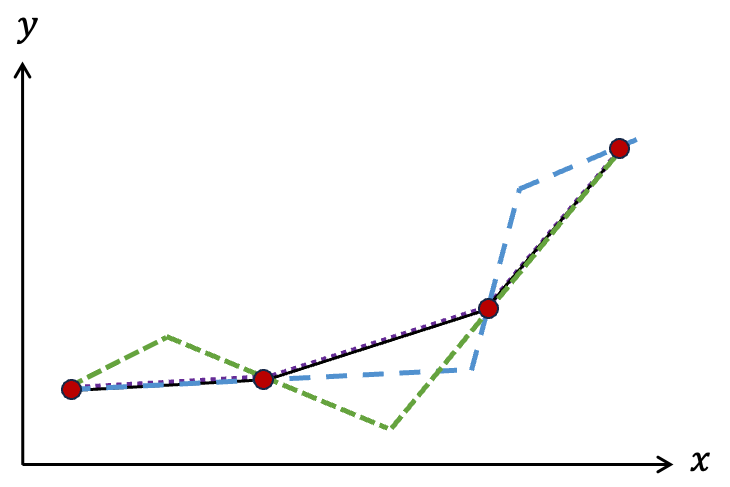}
		\captionof*{figure}{(a) $ d_d < d_d^*, d_h = d_h^* $}
	\end{minipage}%
	\begin{minipage}{.55\textwidth}
		\centering
		\includegraphics[width=.95\linewidth]{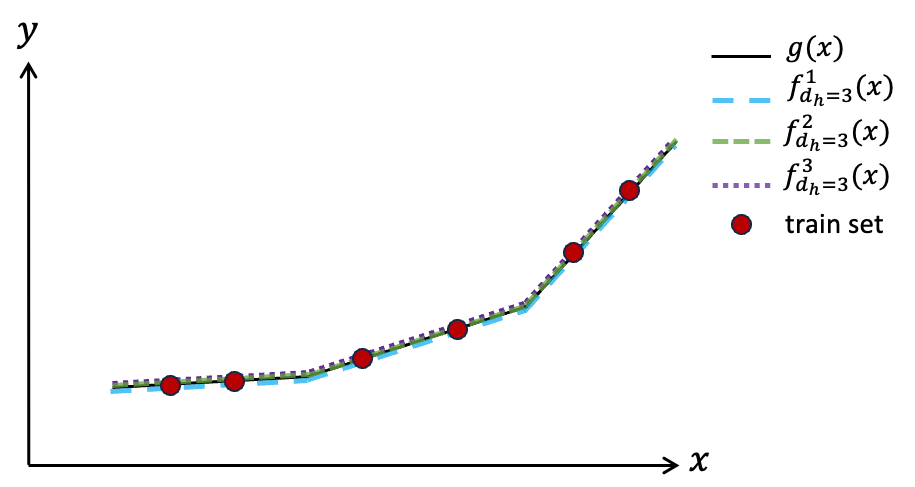}
		\captionof*{figure}{(b) $ d_d = d_d^*, d_h = d_h^* $ \qquad}
	\end{minipage}
    
	\begin{minipage}{.45\textwidth}
		\centering
		\includegraphics[width=.95\linewidth]{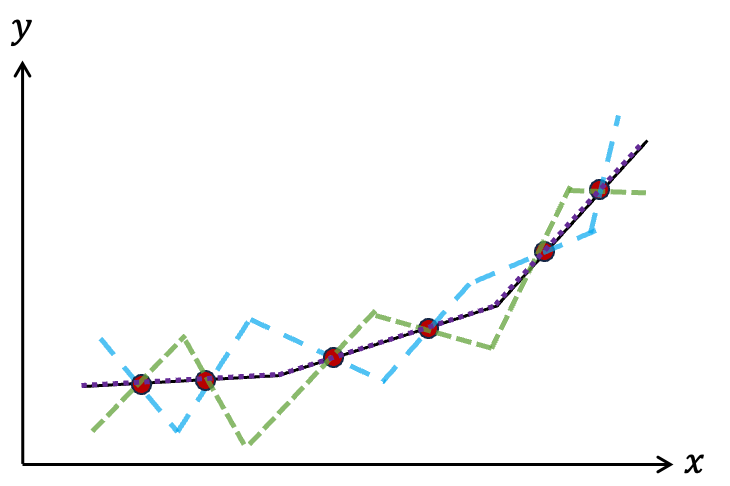}
		\captionof*{figure}{(c) $ d_d = d_d^*, d_h > d_h^* $}
	\end{minipage}%
	\begin{minipage}{.55\textwidth}
		\centering
		\includegraphics[width=.95\linewidth]{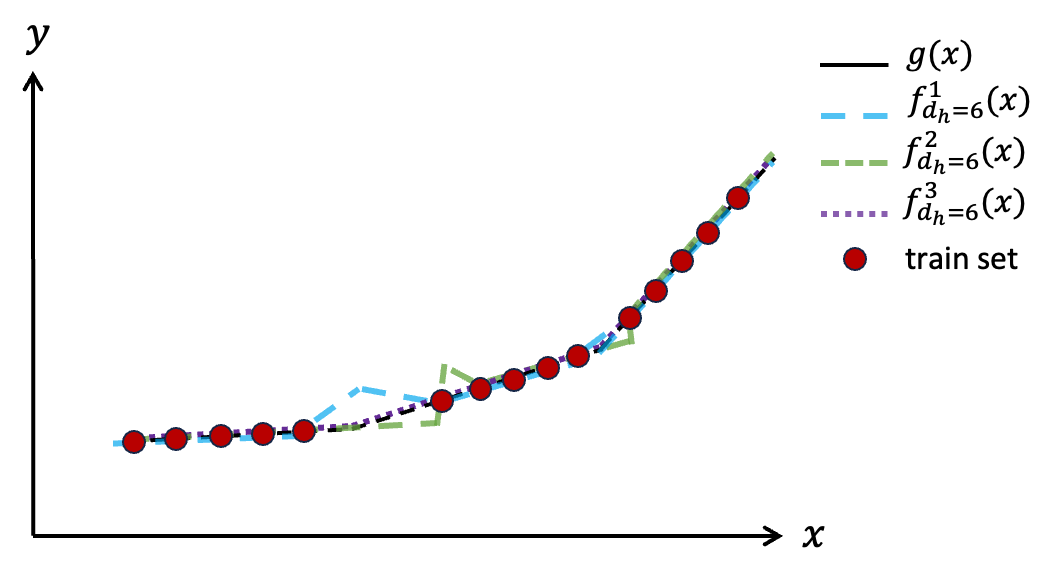}
		\captionof*{figure}{(d) $ d_d > d_d^*, d_h > d_h^* $}
	\end{minipage}
	\caption{The tradeoff between dataset size and model complexity in estimating arbitrary PWL functions}
	\label{fig:armsrace}
\end{figure*}

Figure \ref{fig:armsrace} uses a simple thought experiment to illustrate the trade-off between size of the dataset and the complexity of the model. The target function in the figures, $g(x)$, consists of three linear segments with positive slopes and is estimated by three shallow NN models. Based on the discussion in Section \ref{sec:NNs}, a shallow NN model with $d_h=3$, i.e., $f_{d_h=3}(x)$ is sufficient to represent a PWL function consisting of three linear segments exactly. 

In Figure \ref{fig:armsrace}a, three models $f_{d_h=3}^1(x)$, $f_{d_h=3}^2(x)$ and $f_{d_h=3}^3(x)$ are trained on dataset of 4 data points consisting of the boundary points of each linear segment. The figure illustrates that all three models can reduce the training error to zero by passing through the training points, but only one of them represents the true function $g(x)$. Notice that there are infinitely many different NN models that can zero the training error. 

On the contrary, in Figure \ref{fig:armsrace}b, two data points are sampled for each linear segment. This is enough to limit the model parameter space to a unique set of weights and guarantees that the trained NN model represents $g(x)$ exactly. Notice that simply sampling any 6 data points from $g(x)$ is not enough to guarantee that the model will represent the target function exactly; at least two data points must be sampled from each linear segment. This can be seen from the case where all data points are sampled from the first segment only, in which case the model would not have any information to represent the rest of the function.

It is a common practice in the ML literature to train more complex models than necessary, but this approach has some drawbacks when the target function is deterministic. In Figure \ref{fig:armsrace}c, three models with $d_h=6$ are trained on the same 6 data points as in Figure \ref{fig:armsrace}b. Although each model can represent a PWL function consisting of up to 6 linear segments, they are unable to zero the test error because there are infinitely many ways for the model parameters to be adjusted to minimize the MSE, as illustrated in Figure \ref{fig:armsrace}a. In this case, one needs to increase the number of data points sampled from each segment, here to 5, to reduce the training error, however even this would not guarantee the model to represent the target function.

Other common ML practices used to limit model variance include terminating training if the validation loss no longer reduces or using regularization terms. When the target function is deterministic, early stopping is not a viable option for representing the target function exactly because such a model would not reduce the training error to zero. It can be even better if the validation set examples were added to the training set and overfit the function as essentially the goal is to fit the target function as closely as possible. On the other hand, using regularization with an L1 penalty term can potentially reduce model complexity to the minimum with the proper choice of penalty parameter. In Section \ref{sec:method} we propose an approach that results in a minimum complexity model by design.

% \section{Minimal Size for Exact Estimations}
% \textbf{(How many data points are enough to train a ReLU model to predict a PWL Function)}

\section{A General Framework to Solve QP Using NN}\label{sec:method}
%		In this section, ...
%		\subsection{Modelling QP with Neural Network}
The general form of  mp-QP with linear constraints can be written as,
\begin{subequations}
	\begin{align}
		\textbf{Minimize }  &\  z(\bfm\theta) = \textbf{x}^T\textbf{Q}\textbf{x}+(\textbf{C}+\bfm\theta_c)^T\textbf{x}+\textbf{C}_0,\\
		\textbf{s.t. } &\textbf{A}_e\textbf{x}=\textbf{b}_e + \bfm\theta_c \quad [\lambda],\\
		&\textbf{A}_\mathcal{C}\textbf{x}\leq \textbf{b}_\mathcal{C} + \bfm\theta_\mathcal{C} \quad [\mu],
	\end{align}
	\label{eq:ineq}
\end{subequations} 
where $\textbf{x}, \textbf{C}, \textbf{C}_0 \in \mathbb{R}^{n}, \textbf{b}_e\in \mathbb{R}^{m_1}, \textbf{b}_\mathcal{C}\in \mathbb{R}^{m_2}$, $\textbf{Q} \in \mathbb{R}^{n\times n}$ is a semi-positive definite matrix, and $\textbf{A}_e\in \mathbb{R}^{m_1\times n}$, $\textbf{A}_\mathcal{C}\in \mathbb{R}^{m_2\times n}$. The problem is parametrized by $\bfm\theta = [\bfm\theta_c,\bfm\theta_e,\bfm\theta_\mathcal{C}]$ for $\bfm\theta_c\in \mathbb{R}^{n},\bfm\theta_e\in \mathbb{R}^{m_1},\bfm\theta_\mathcal{C}\in \mathbb{R}^{m_2}$ and $\bfm\theta \in \Theta_f$, where $\Theta_f$ is the set of all feasible parameters. The Lagrangian function of the above problem can be written as,
\begin{equation}
	L(\textbf{x},\boldsymbol{\lambda},\boldsymbol{\mu}) = \textbf{x}^T\textbf{Q}\textbf{x}+(\textbf{C}+\bfm\theta_c)^T\textbf{x}+\textbf{C}_0 +\boldsymbol{\lambda}^T( \textbf{b}_e+\bfm\theta_e-\textbf{A}_e\textbf{x})+\boldsymbol{\mu}^T( \textbf{b}_\mathcal{C}+\bfm\theta_\mathcal{C}-\textbf{A}_\mathcal{C}\textbf{x}),
	\label{eq:lagr}
\end{equation}
where $\boldsymbol{\lambda} \in \mathbb{R}^{m_1},\boldsymbol{\mu}\in \mathbb{R}^{m_2}$ are the dual variables associated with equality and inequality constraints, respectively. 

If the solution to the subproblem (\ref{eq:ineq}a)-(\ref{eq:ineq}b)  satisfies (\ref{eq:ineq}c), the corresponding shadow prices must be zero at optimality, i.e., $\boldsymbol{\mu}^*=0$. However, when the solution violates (\ref{eq:ineq}c), $\bfm\mu^*$ can be obtained by adding one or more violated constraints iteratively to reach the optimal solution. In the following sections, the general cases of optimization problems with equality constraints (\ref{eq:ineq}b) and inequality constraints (\ref{eq:ineq}c) are considered, and methods for estimating the dual variables are proposed. Finally, a NN based solution method to such optimization problems is proposed. 
%    Therefore, there is always a set of binding inequality constraints, empty or non-empty, that can be turned to equality constraint. Therefore, if the set of binding constraints could be predicted, the problem could be reduced to an equality constrained problem, whose solution forms a linear function. 

\subsection{Solution to Equality Constrained mp-QP}	\label{sec:model_eq}
For the case when constraint (\ref{eq:ineq}c) is not binding, $ \boldsymbol{\mu}^*=0 $ and the Lagrangian function (\ref{eq:lagr}) can be written as follows.
%	    To solve the optimization problem (\ref{eq:ineq}a)-(\ref{eq:ineq}b),  first the Lagrangian function (\ref{eq:lagr_eq}) of the problem is constructed, whose derivatives yield the optimal solution to the problem.
%	    
%	    In case when the inequality constraints are not binding, constraint (\ref{eq:ineq}c) can be ignored as $\boldsymbol{\mu}^*=0$, and therefore (\ref{eq:lagr}) can be written as below:
\begin{equation}
	L(\textbf{x},\boldsymbol{\lambda},\boldsymbol{\mu}) = \textbf{x}^T\textbf{Q}\textbf{x}+(\textbf{C}+\bfm\theta_c)^T\textbf{x}+\textbf{C}_0 +\boldsymbol{\lambda}^T( \textbf{b}_e+\bfm\theta_e-\textbf{A}_e\textbf{x}).
	\label{eq:lagr_eq}
\end{equation}
Using the fact that the gradient of the Lagrangian function (\ref{eq:lagr_eq}) with respect to $\textbf{x},\boldsymbol{\lambda}$ must be zero at optimality, equation (\ref{eq:linearsolver}) defines a function that produces optimal solutions to any given right hand side vector, $ \textbf{b}_e $ (see  \ref{app:linear_derivation} for derivation).
%From the above equation \ref{eq:linearinverse},  the function $ g(\textbf{C},\textbf{b}_1) $ can be obtained as
\begin{align}
	g(0;\bfm\theta) = \textbf{J}^{-1}
	\begin{bmatrix}
		-\textbf{C}-\bfm\theta_c \\ -\textbf{b}_e-\bfm\theta_e
	\end{bmatrix},
	\quad \text{where} \quad \textbf{J} = \begin{bmatrix}
		2\textbf{Q} & -\textbf{A}_e^T\\
		-\textbf{A}_e & 0
	\end{bmatrix},
	\label{eq:linearsolver}
\end{align}
represents the coefficient matrix. Equation (\ref{eq:linearsolver}) is a function yielding the optimal solutions for (\ref{eq:ineq}a)-(\ref{eq:ineq}b) as $g(0;\bfm\theta) = [\textbf{x}^*,\boldsymbol{\lambda}^*]$, which we refer to as the \textit{solution function}. Here the 0 in the $g(\cdot)$ function is the value of the shadow price, $\bfm\mu^*=0$, which will be used to generalize to inequality constrained problems in the next section.
% Based on this observation, the following theorem presents how a NN, trained with optimal solutions, represents the system parameters in the hidden layers.

%		and provide insight into our approach t obtain model parameters without training but fixing them to the weights derived from the optimization problem.

% To find the optimal solution $\textbf{x}^*$ corresponding to each $\textbf{b}_1$, one can derive $g^{-1}(\textbf{b}, \textbf{e})$ via inverting the above Jacobian matrix, $\textbf{J}$. Notice that $g^{-1}(\textbf{b}, \textbf{e})$ is a linear function and, based on Theorem \ref{th:linear_parameters}, this function can be estimated by a NN with an identity activation where the product of optimal parameter weights equals the inverse of the Jacobian matrix. Optimization problems with quadratic objective functions are a special case where the partial derivatives form a linear system of equations. In the next section, we generalize our approach to problems with nonlinear partial derivatives.

\subsection{Adding Inequality Constraints}
As discussed earlier, when inequality constraints are binding, the active constraints can be considered as equality constraints with the optimal dual variables, $\bfm\mu^*>0$. Then, the solution can be found by expanding the Lagrangian function and finding its critical points.

Assume that for any system parameter, $\bfm \theta$, the set of indices of binding constraints $\mathcal{B}\subseteq \{1,\dots,m_2\}$ is known. Then, the optimal solution to the problem can be found by solving
\begin{subequations}
	\begin{align}
		\frac{\partial L}{\partial \textbf{x}} &= 2\textbf{Q}\textbf{x}+\textbf{C}+\bfm\theta_c + \boldsymbol{\lambda}^T\textbf{A}_e+\boldsymbol{\mu}^T\textbf{A}_\mathcal{C} = 0\\
		\frac{\partial L}{\partial \boldsymbol \lambda} &= \textbf{b}_e+\bfm\theta_c-\textbf{A}_e\textbf{x} =0,\\
		\frac{\partial L}{\partial \boldsymbol \mu} &= \textbf{b}_{\mathcal{B}}+\bfm\theta_{\mathcal{B}}-\textbf{A}_{\mathcal{B}}\textbf{x} =0,
	\end{align}
	\label{eq:binding_der}
\end{subequations}
where $\textbf{A}_\mathcal{B} = \{\textbf{A}_\mathcal{C}\}_{i\in\mathcal{B}}, \bfm\theta_\mathcal{B} = \{\bfm\theta_\mathcal{C}\}_{i\in\mathcal{B}}, \textbf{b}_\mathcal{B}= \{\textbf{b}_\mathcal{C}\}_{i\in\mathcal{B}}$ are matrices consisting of binding constraints, and optimal shadow prices $\mu_i^* > 0 \  \forall i \in \mathcal{B} $ and $\mu_i^* = 0 \  \forall i\notin \mathcal{B} $.

Now, if the function $\bfm\mu(\bfm\theta)$ that generates $\bfm\mu^*$ is known, the remaining optimal solutions can be obtained by modifying (\ref{eq:lagr}) and writing the derivatives as
\begin{subequations}
	\begin{align}
		\frac{\partial L}{\partial \textbf{x}} &= 2\textbf{Q}\textbf{x}+\textbf{C}+\bfm\theta_c + \boldsymbol{\lambda}^T\textbf{A}_e+\boldsymbol{\mu}^T(\bfm\theta)\textbf{A}_\mathcal{C} = 0,\\
		\frac{\partial L}{\partial \boldsymbol \lambda} &= \textbf{b}_e+\bfm\theta_c-\textbf{A}_e\textbf{x} =0.
	\end{align}
\end{subequations}
The above mapping is the same as the equality constrained solution function in \ref{eq:linearsolver}. Therefore, if $\bfm\mu(\bfm\theta) = \boldsymbol{\mu}^*$ is given, the optimal solution can be obtained using the solution function,
\begin{equation}
	g(\bfm\mu^*;\bfm\theta) = [\textbf{x}^*,\boldsymbol{\lambda}^*]^T.
	\label{eq:piecewisesolver}
\end{equation}
As the solution function is linear, the only nonlinearity can stem from $\boldsymbol{\mu}^*$. The next section discusses the piecewise linear nature of $\boldsymbol{\mu}^*$ and how it can be estimated exactly with a NN.

\subsection{Predicting Active Constraints and Shadow Prices}
As shown in (\ref{eq:piecewisesolver}), estimating the solution of the optimization problem (\ref{eq:ineq}a)-(\ref{eq:ineq}b) requires knowing the value of $ \mu^* $. Therefore, in this section, a methodology to produce exact predictions of $\boldsymbol{\mu}^*$ and hence the active constraints is presented. While $ g(\cdot) $ in eq. (\ref{eq:piecewisesolver}) is a linear function, the addition of $ \mu^* $ changes the slope of the solutions with respect to $\textbf{b}_e$ and converts the function to PWL due to the piecewise linear nature of $ \mu(\bfm\theta) $.

\begin{theorem}\label{th:pwa}
	For the mp-QP problem \ref{eq:ineq}, $\Theta_f \subseteq \Theta$ is a convex set, the primal solution function \textbf{x}$(\bfm\theta): \Theta_f \rightarrow \mathbb{R}^n$ is continuous and piecewise affine. Also the optimal objective function $\textbf{z}(\bfm\theta):\Theta_f\rightarrow \mathbb{R}$ is continuous and piecewise quadratic.
\end{theorem}
\textit{Proof:} The reader can find the proof of the theorem in \cite{pistikopoulos2020multi}.

\begin{lemm}
	Consider the mp-QP problem \ref{eq:ineq}. The shadow prices of inequality constraints, $\bfm\mu^*:\Theta_f \rightarrow \mathbb{R}^{m_2}$ are piecewise affine.
\end{lemm}
\textit{Proof:} The mapping $g(\cdot)$ in eq. \ref{eq:piecewisesolver} is an affine mapping from $\bfm\mu^*=\mu(\bfm\theta)$, and as shown in Theorem \ref{th:pwa}, the output of $\textbf{x}(\bfm\theta)$ is a piecewise affine function of $\bfm\theta$. As $g(\cdot)$ can be defined as a composite function, i.e., $(g\circ \mu)(\bfm\theta)$, the only input of the function has to be piecewise affine.

To illustrate, consider an optimization problem in the form defined in (\ref{eq:ineq}), with 5 decision variables, $x_i$, three equality constraints, and three upper and lower bound constraints; i.e., $ x_i^{-} \leq x_i \leq x_i^{+} $ for $i \in [0,1,2]$. As will be shown later, this is our smallest 6-bus DC-OPF application case defined in eq. \ref{eq:DC-OPF}. Here, the problem is solved using a Gurobi solver (v10.0.3) for $\theta_1\in [0,0.1,0.2,\dots,5]$, controlling the RHS of equality constraints $ \textbf{b}_e + \bfm\theta_c= [\theta_1,0.001,0.001] $. In Figure \ref{fig:mu}, the optimal dual variables are plotted against different values of $ \theta_1 $. $\mu_i$ for $i \in [0,1,2]$ and $\mu_i$ for $i \in [3,4,5]$ are shadow prices reflecting upper and lower power generator limits, respectively, in the example. The slope changes in $\bfm\mu^*$ exhibit a piecewise affine pattern, with different slopes in different critical regions (\textbf{CR}), corresponding to different combinations of active constraints. In $\textbf{CR}_0$, two lower limit constraints are active, with corresponding $\mu_i>0$. As $\theta_1$ increases, generator lower limits stop binding and upper limits start to bind.

While $\boldsymbol{\mu}^*$ as a function of $ \bfm\theta $ is piecewise linear and can be modelled with a generic NN, training a model on a dataset consisting of $(\bfm\theta, \boldsymbol{\mu}^*)$ pairs  does not guarantee optimality of the predictions for test data sampled outside this distribution. To overcome this difficulty, our method directly injects some of the information that is general to all $\bfm\theta\in \Theta_f$, into the first layer of the model. Specifically, all the slopes that the $\mu$-NN will have for each potential combination of binding constraints are calculated and fixed as the first layer weights prior to the training, then the rest of the model parameters, biases and second layer weights, are estimated via training. Figure \ref{fig:muPredictorNN} illustrates the architecture of this model where $\textbf{W}^0$ indicates the injected weights for the first layer, obtained from $\nabla \boldsymbol{\mu}^*_\mathcal{B}$ as follows:

\begin{figure}
	\centering
	\begin{minipage}{.5\textwidth}
		\centering
		\includegraphics[width=.8\linewidth]{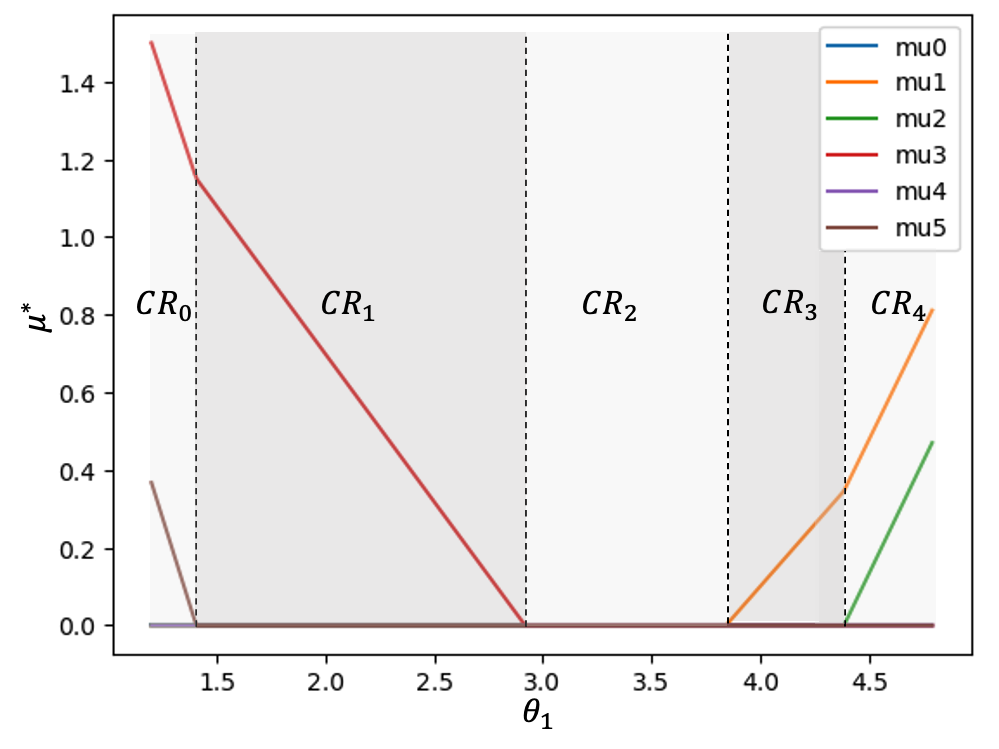}
		\caption{Sensitivity of $\boldsymbol{\mu}$ with respect to $ 
			\theta_1$}
		\label{fig:mu}
	\end{minipage}%
	\begin{minipage}{.5\textwidth}
		\centering
		\vspace{25pt}
		\includegraphics[width=.9\linewidth]{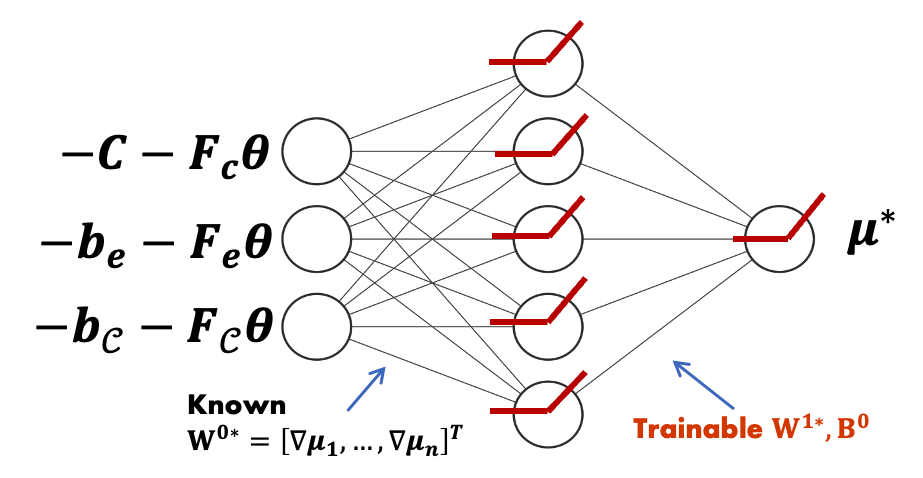}
		\caption{NN Model to Predict $\boldsymbol{\mu}^*$}
		\label{fig:muPredictorNN}
	\end{minipage}
\end{figure}

\textit{Calculation of} $\bfm{\mu}^*_{\mathcal{B}_i}$: To obtain the slopes of $\boldsymbol{\mu}^*$ with respect to the changes in $\bfm \theta$ for binding constraints, $\nabla \boldsymbol{\mu}^*_{\mathcal{B}_i}$ we expand the \textbf{J} matrix defined in eq. \ref{eq:linearsolver} with the rows of the coefficient matrix of inequality constraints corresponding to the binding constraints, $\textbf{A}_{\mathcal{B}_i}$, and then invert it.
\begin{align}
	\begin{bmatrix}
		\textbf{x}^*\\ \boldsymbol{\lambda}^* \\ \boldsymbol{\mu}^*_{\mathcal{B}_i}
	\end{bmatrix} =
	\textbf{J}_{\mathcal{B}_i}^{-1}
	\begin{bmatrix}
		-\textbf{C} - \bfm\theta_c \\ -\textbf{b}_e-\bfm\theta_c \\ -\textbf{b}_{\mathcal{B}_i}-\bfm\theta_{\mathcal{B}_i}
	\end{bmatrix}
	\begin{bmatrix}
		\nabla \textbf{x}^* \\ \nabla \boldsymbol{\lambda}^* \\ \nabla \boldsymbol{\mu}^*_{\mathcal{B}_i}
	\end{bmatrix}
	\begin{bmatrix}
		-\textbf{C} - \bfm\theta_c \\ -\textbf{b}_e-\bfm\theta_c \\ -\textbf{b}_{\mathcal{B}_i}-\bfm\theta_{\mathcal{B}_i}
	\end{bmatrix},
	\label{eq:linearinverse_ineq}
\end{align}
where, 
\begin{align}
	\textbf{J}_\mathcal{B} = \begin{bmatrix}
		2\textbf{Q} & -\textbf{A}_e^T & -\textbf{A}_{\mathcal{B}_i}^T\\
		-\textbf{A}_e & 0 & 0\\
		-\textbf{A}_{\mathcal{B}_i} & 0 & 0
	\end{bmatrix}.
\end{align}

Notice that the last row of the above inverse Jacobian matrix, $\textbf{J}^{-1}_{\mathcal{B}_i}$, is the gradient vector $\nabla \boldsymbol{\mu}^*_{\mathcal{B}_i}$ that contains the coefficients of the linear segment of the target PWL function when a certain set of constraints, $\mathcal{B}_i$, are binding. In other words, the vector is the set of coefficients that maps $[\textbf{C},\textbf{b}_e, \textbf{b}_{\mathcal{B}_i}]$ and $\bfm \theta $ to $\boldsymbol{\mu}^*_{\mathcal{B}_i}$.

%	To express mathematically, let $\mathcal{J}_\mu^{-1}$ be the rows of the inverse Jacobian corresponding to $\boldsymbol{\mu}$,
%	 	\begin{align}
	%		\mathcal{J}_\mu^{-1} \begin{bmatrix}
		%			-\textbf{C}\\ \textbf{b}_1 \\ \textbf{b}_2^\mathcal{B}
		%		\end{bmatrix} = \boldsymbol{\mu}^*_\mathcal{B} \Rightarrow \left[-\frac{\partial \boldsymbol{\mu}}{\partial \textbf{C}},\frac{\partial \boldsymbol{\mu}_\mathcal{B}}{\partial \textbf{b}_1},\frac{\partial \boldsymbol{\mu}}{\partial \textbf{b}_2^\mathcal{B}}\right] = \mathcal{J}_\mu^{-1}.
	%	\end{align}
After the first layer of the NN is fixed with all potential derivatives, the training is carried out using $(\bfm\theta,\bfm\mu^*)$ pairs to find when each slope is activated. However, the number of potential slopes increases exponentially with the number of inequality constraints, which increases the computational burden of the pre-training step as matrix inversion is required to calculate each slope, and also increases the number of parameters to estimate in the training step. The next section outlines a method to identify critical regions of the problem and describes a data augmentation procedure that reduces the usage of solvers.

\subsection{Identification of Critical Regions and Data Generation}
% The proposed method can create a NN model that can produce error-free predictions for the shadow prices, but there are practical challenges in training the model, because it is necessary to find the slopes of all cases of binding constraints to fix the first layer weights. Specifically, for a system with $m_2$ inequality constraints, the number of possible set of binding constraints is $2^{m_2}.$

\begin{figure}
	\centering
	\begin{minipage}{.55\textwidth}
		\centering
		\vspace{5pt}
		\includegraphics[width=.8\linewidth]{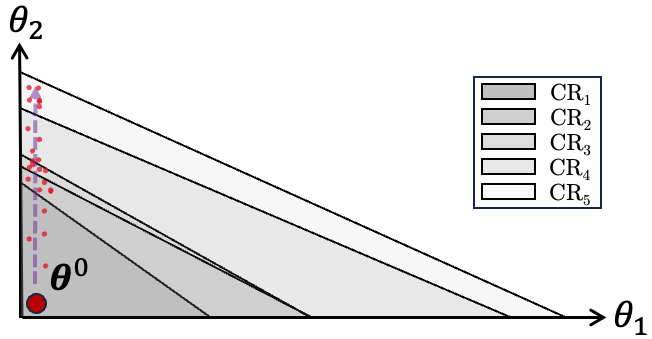}
		\caption*{(a) PSNN}
	\end{minipage}%
	\begin{minipage}{.55\textwidth}
		\centering
		\vspace{2pt}
		\includegraphics[width=.8\linewidth]{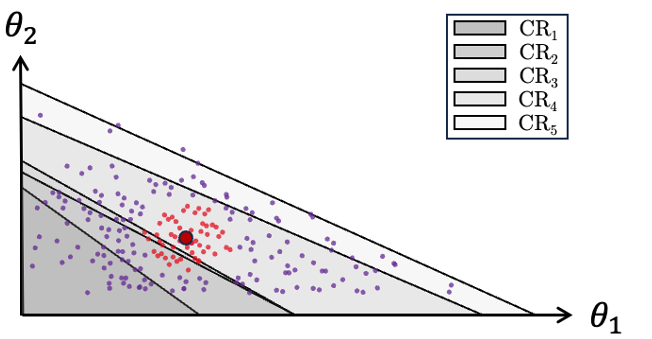}
		\caption*{(b) DNN}
	\end{minipage}
	\caption{Illustration of critical region discovery and training set sampling for PSNN and DNN models}
	\label{fig:discovery}
\end{figure}

As explained in the Introduction, the problem characteristics in many mp-QP applications imply that most constraint combinations can never occur, so some method of filtering out the unused combinations is needed. For example, in the above 6-bus DC-OPF case with 6 inequality constraints, the set of all constraint combinations is $\{\{0\}, \{1\},\{2\}, \{0,1\}, \dots, \{0,1,2,3,4,5\}\}$ with $2^6=64$ elements, yet 59 of these combinations never occur. In the power demand applications investigated in the present study, we use a commercial solver to identify the active constraint combinations by solving the problems over a range of input parameters.

begins with an input parameter, \(\theta^0\), and solves the problem using a commercial solver at \(\theta^0\) to identify the active set of constraints at the initial solution, \(\mathcal{B}_0\). 

Algorithm \ref{alg:discovery} details this discovery step. The algorithm first generates data starting from an initial parameter, $\bfm\theta^0$ and choosing an arbitrary dimension, $k$, and incrementally increasing one dimension from zero to the highest possible feasible value. Then, after solving for the first input demand, the set of potential binding constraints, \(\mathcal{B}_0\), is identified by checking whether the associated $\mu_i$ is greater than zero for the inputs. Using this active set, it expands the coefficient matrix and obtains slope coefficients $\nabla\bfm\mu_{\mathcal{B}_i}$ and the matrix \(J^{-1}_{\mathcal{B}_0}\) . This matrix is then used to calculate solutions for the input parameters and to check whether the obtained solutions satisfy the KKT conditions. If the conditions are not satisfied, a new critical region is identified, and the problem is solved again by the solver to determine the new active set, \(\mathcal{B}_{i+1}\). The remaining input parameters are then solved using the updated solution function, \(J^{-1}_{\mathcal{B}_{i+1}}\). This procedure is repeated until all critical regions have been discovered for all input parameters. Whenever a critical region is identified, the algorithm also provides insight into the boundaries of critical regions, $[\textbf{CR}_i^-,\textbf{CR}_i^+]$ corresponding to $\mathcal{B}_i$ that are active at solution. This information is used in Algorithm \ref{alg:populate} to populate data without using solvers, by sampling random inputs from $[\textbf{CR}_i^-,\textbf{CR}_i^+]$ for all $i$ and calculating the corresponding optimal solutions using the slope corresponding to each critical region.

Finally, Algorithm \ref{alg:train} details the training of $\mu$-NN. The first layer of the model is fixed to $\textbf{W}^0$ and and the rest of the model parameters, $\textbf{W}^1,\textbf{B}^0$, are randomly initiated. The model is trained with 1000 datapoints that are constructed as detailed in Algorithm \ref{alg:populate}.

% \sloppy
% The training procedure can be summarized as follows (see  \ref{app:training} for further details):
% \begin{enumerate}
% 	\item Generate RHS and solve with a solver: (i) Initiate the RHS with an initial point, (ii) Choose a $\theta_k$ and increase its value from 0 to max.
% 	\item Solve the problem and find $\textbf{CR}_i$ and $\nabla \boldsymbol{\mu}^*_{\mathcal{B}_i}$:
% 	(i) Solve the problem for each generated data using a commercial solver, (ii) Find the sets of active constraints, i.e., $\mathcal{B}_i = \{\ j \ | \ \mu_j>0 \ \}$, (iii) Expand the coefficient matrix and calculate $\textbf{J}^{-1}_{\mathcal{B}_i}$ for all binding constraints. Obtain $\nabla \boldsymbol{\mu}^*_{\mathcal{B}_i}$.
% 	\item Populate data: Randomly generate $d_d$ number of RHS vectors using $Uniform(\textbf{CR}_i^-, \textbf{CR}_i^+)$ and find optimal solutions following \ref{eq:linearinverse_ineq}. Check if KKT is satisfied for each data and regenerate data if needed.
% 	\item Train $\mu$-NN: (i) Initiate $\mu$-predictor NN with random $\textbf{W}^1,\textbf{B}^0$, fix $\textbf{W}^0 = [\boldsymbol{\mu}^*_{\mathcal{B}_1},\dots,\boldsymbol{\mu}^*_{\mathcal{B}_n}]^T$ and set with random weights, (ii) Train the $\mu$-NN with 1000 datapoints populated using the procedure described above, (iii) Define solver NN by setting first part to $\mu$-NN and second part to $\textbf{J}^{-1}$.
% \end{enumerate}

\begin{algorithm} \label{alg:discovery}
    \caption{Data Generation and Discovery}
    \label{alg:discovery}
    \begin{algorithmic}
    \State \textbf{Input:} Initial point $\bfm\theta^0$, increment size $\alpha$, upper limit $\theta^+$, KKT violation threshold $tol$
    \State \textbf{Output:} ($\textbf{CR}_i^-,\textbf{CR}_i^+$) and $\mathcal{B}_i$ for all $i$
    \State Initiate an empty $\mathcal{D}$
    \State $\bfm\theta \gets \bfm\theta^0$ 
    \For{$d$ in 0 to $\theta^+$ step $\alpha$}
    \State Append $\bfm\theta_k$ to $\mathcal{D}$ for $\bfm\theta_k = d$
    \EndFor
    \State Solve problem for the initial point using a solver 
    \State Find active constraints: $\mathcal{B}_0 = \{\ j \ | \ \mu_j > 0 \ \}$
    \State Expand the coefficient matrix, calculate $\textbf{J}^{-1}_{\mathcal{B}_0}$ and $\nabla\bfm\mu_{\mathcal{B}_i}^*$
    \State Initiate $i=0$, $\textbf{CR}^-=[\ ]$, $\textbf{CR}^+=[\ ]$, $\textbf{W}^0=[\ ]$, 
    \State Set $\textbf{CR}^-[i]=\bfm\theta$, $\textbf{CR}^+[i]=\bfm\theta$ and $\textbf{W}^0[i]= \nabla\bfm\mu_{\mathcal{B}_i}$, 
    \For{$\bfm\theta$ in $\mathcal{D}$}
    \State Calculate $[\textbf{x}^*,\bfm\lambda^*,\bfm\mu^*_{\mathcal{B}_i}]^T$ for $\bfm\theta$ by substituting $\textbf{J}^{-1}_{\mathcal{B}_i}$ in (13)
      \State Check if $[\textbf{x}^*,\bfm\lambda^*,\bfm\mu^*_{\mathcal{B}_i}]^T$ satisfies all KKT conditions in Table 1 with less than $tol$ error
      \If{KKT is not satisfied, a new critical region is found,}
      \State Set $i=i+1$
      \State (i) Solve the problem for $\bfm\theta$ using a solver
      \State (ii) Find active constraints: $\mathcal{B}_i = \{\ j \ | \ \mu_j > 0 \ \}$
      \State (iii) Expand coefficient matrix to calculate $\nabla\bfm\mu^*_{\mathcal{B}_i}$ and append to $\textbf{W}^0$
      \EndIf
      \State Update $\textbf{CR}^+[i]=\bfm\theta$
      % \If{$\mathcal{B}_i\notin \textbf{W}^0$}
      % \State Expand coefficient matrix to calculate $\nabla\bfm\mu_{\mathcal{B}_i}$ and append to $\textbf{W}^0$
      % \State Set $\textbf{CR}_i^+=\bfm\theta$
      % \State Set $i=i+1$
      % \EndIf
    \EndFor\\
    \Return $\textbf{W}^0, \textbf{CR}^-,\textbf{CR}^+$
    \end{algorithmic}
\end{algorithm}

    \begin{algorithm}
    \caption{Training Data Generation}
    \label{alg:populate}
    \begin{algorithmic}
    \State \textbf{Input:} $(\textbf{CR}^-_i,\textbf{CR}^+_i), \nabla \boldsymbol{\mu}^*_{\mathcal{B}_i}$ for all $i \in [1,\dots,n]$ and $d_{d}$
    \State \textbf{Output:} $(\bfm\theta, \bfm\mu^*)$ pairs.
          \For{$i$ in [1,\dots,n]}
          \State (i) Randomly generate $d_d$ number of $\bfm\theta$ vectors by sampling $d\sim $ Uniform($\textbf{CR}_i^-$, $\textbf{CR}_i^+$) and setting kth dimension of $\bfm\theta$ with each $d$
          \State (ii) Find optimal solutions as $\bfm\mu^* = \nabla \boldsymbol{\mu}^*_{\mathcal{B}_i}(-\textbf{B}-\bfm\theta)$
          \State (iii) Check if KKT is satisfied for each data
          \If {KKT not satisfied}
            \State For each unsatisfied example, repeat (i)-(iii)
        \EndIf
      \EndFor
    \end{algorithmic}
\end{algorithm}

\begin{algorithm}
\caption{Train $\mu$-NN}
\label{alg:train}
    \begin{algorithmic}
    \State \textbf{Input:} $\textbf{W}^0 = [\nabla\bfm\mu_{\mathcal{B}_1},\dots,\nabla\bfm\mu_{\mathcal{B}_n}]$, maximum epoch $M$, learning rate $\eta$.
    \State \textbf{Output:} Labeled dataset, $\textbf{CR}_i$ and $\mathcal{B}_i$ $\forall i$
    \State Initiate $\mu$-NN with random weights, $\textbf{W}^1$, $\textbf{B}^0$
    \State Fix the first layer weights to $\textbf{W}^0 = [\boldsymbol{\mu}^*_{\mathcal{B}_1},\dots,\boldsymbol{\mu}^*_{\mathcal{B}_n}]^T$
    \State Train the model with $(\bfm\theta,\bfm\mu)$ pairs with $\eta$ learning rate until either (i) $MSE_{val}<tol$ or (ii) $epoch==M$
    \end{algorithmic}
  \end{algorithm}

Figure \ref{fig:discovery} compares the search pattern and dataset construction for PSNN and traditionally trained DNN models. In Figure \ref{fig:discovery}a the arrow represents the search pattern where the data points are generated equidistantly and red points represent the minimal number of data points sampled  for every critical region. As seen in the figure, moving along one axis is sufficient to visit all critical regions and identify the bounds of these regions in our examples. On the other hand, the dataset for DNN training is sampled by perturbing an initial (average) parameter set. It is shown in the next section that when PSNN is trained on data points sampled along one axis, it is capable of predicting optimal solutions for all feasible region, whereas DNN model does not generalize outside this domain.

Note that different problems have different topologies and this search pattern would not be sufficient to discover all regions. Nevertheless, the performance of our model on the test cases is important to show the generalization power outside the training distribution thanks to the analytical derivation of the slopes. This point is addressed further in the limitations section.

\section{Numerical Results}\label{sec:dcsolver}
% This section will detail the DC-OPF problem, derive the equations to be modelled following the procedure described above and discuss the IEEE test cases used to compare the performance of our model against alternative approaches. 
% %		Gurobi solver results and a NN model trained with a dataset constructed by solving the problem repeatedly with a solver. Next section details DC-OPF problem, followed by the results on different test systems.

% \subsection{DC-OPF problem}\label{sec:dc}
%	Optimal Power Flow is a nonlinear and nonconvex problem to find the optimal generator dispatches that minimize the total cost of power production and satisfies the power demand by respecting the grid system constraints such as ramp-up and ramp-down constraints, line limits and other operational constraints. The problem is challenging due to its nonlinear and nonconvex nature and is typically solved with computationally costly solvers. To be able to operate the system, the power flow constraints are linearized with some assumptions, which results in a QP with linear constraints.
Optimal Power Flow (OPF) is a central problem in electrical power system management \citep{carpentier1979optimal}. The problem is to find the optimal
dispatch per generator to minimize an objective function (e.g., cost of operation) while satisfying energy demands and grid system constraints at a given point in time. The formulation of the AC-OPF (alternating current) problem is nonlinear and nonconvex, so challenges in finding optimal and feasible solutions led researchers to use the simplified DC-OPF (direct current) formulation which has linearized power flow constraints under certain assumptions. 
% The motivation behind the linear relaxation is that the difference between voltage angles at transmission level is very small, i.e. $\theta_i-\theta_j \approx 0$. Hence it can be assumed that $\sin(\theta_i-\theta_j) \approx \theta_i -\theta_j$ and $\cos(\theta_i -\theta_j) \approx 1$. 
The DC-OPF problem can be defined as:
%	DC-OPF is a QP with linear constraints that aims to find the optimal generator dispatches that satisfies the power demand 
\begin{subequations}
	\begin{align}
		\underset{\textbf{P}_g,\boldsymbol{\delta}}{\textbf{Minimize: }}  & \textbf{P}_g^T\textbf{Q}\textbf{P}_g+\textbf{C}^T\textbf{P}_g+\textbf{C}_0,\\
		\textbf{s.t. } &\textbf{P}_d +\bfm\theta_c - \textbf{B}\boldsymbol{\delta}-\textbf{P}_g=0 \quad [\boldsymbol{\lambda}],\\
		&\textbf{P}_g-\textbf{P}_g^+\leq 0 \qquad \quad\qquad[\boldsymbol{\mu}^+], \\
		&\textbf{P}_g^--\textbf{P}_g\leq 0 \qquad\quad\qquad[\boldsymbol{\mu}^-],
	\end{align}
	\label{eq:DC-OPF}
\end{subequations}
where $\textbf{C}, \textbf{C}_0 \in \mathbb{R}^{n_g}$ are cost coefficient vectors, $\textbf{Q} \in \mathbb{R}^{n_g\times n_g }$ is the quadratic cost coefficient matrix, $\textbf{P}_g^{-}, \textbf{P}_g^{+}\in \mathbb{R}^{n_g}$ are lower and upper generator limits, and $n_g, n_b$ are generator buses and the total number of buses, respectively. The primal variables are the production output of dispatchable generators, $\textbf{P}_g  \in \mathbb{R}^{n_g}$ and voltage magnitudes of each bus, $\boldsymbol{\delta} \in \mathbb{R}^{n_b}$. Following the steps in Section \ref{sec:method}, our model is derived using the derivatives of the Lagrangian function in eq. \ref{eq:binding_der} for the DC-OPF problem. The reader can refer to \ref{app:dc_derivation} for the derivation of the system of equations. After the derivations, the right hand side vectors translate into the demand vector, $\textbf{b}_e=\textbf{P}^{\text{lim}}=\textbf{P}_d$ for equality constraints and the generator limits, $\textbf{b}_\mathcal{C}^T=[\textbf{P}_g^{+},\textbf{P}_g^{-}]$ for the inequality constraints, and the primal variables are $\textbf{x}^T=[\textbf{P}_g,\boldsymbol{\delta}]$.

\begin{figure}
	\centering
	
	\begin{minipage}{.45\textwidth}
		\centering
		\vspace{5pt}
		\includegraphics[width=1\linewidth]{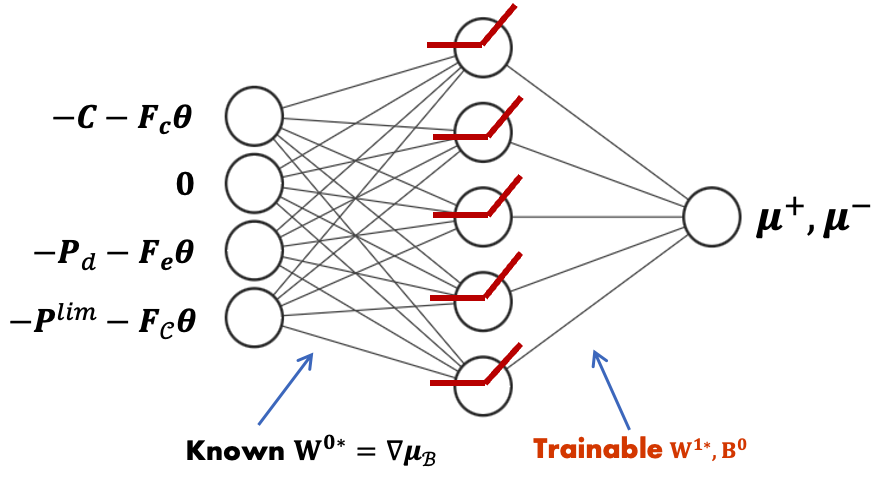}
		\caption{Subnetwork predicting $\boldsymbol{\mu}^*$}
		\label{fig:mupred}
	\end{minipage}%
	\begin{minipage}{.45\textwidth}
		\centering
		\vspace{10pt}
		\includegraphics[width=1\linewidth]{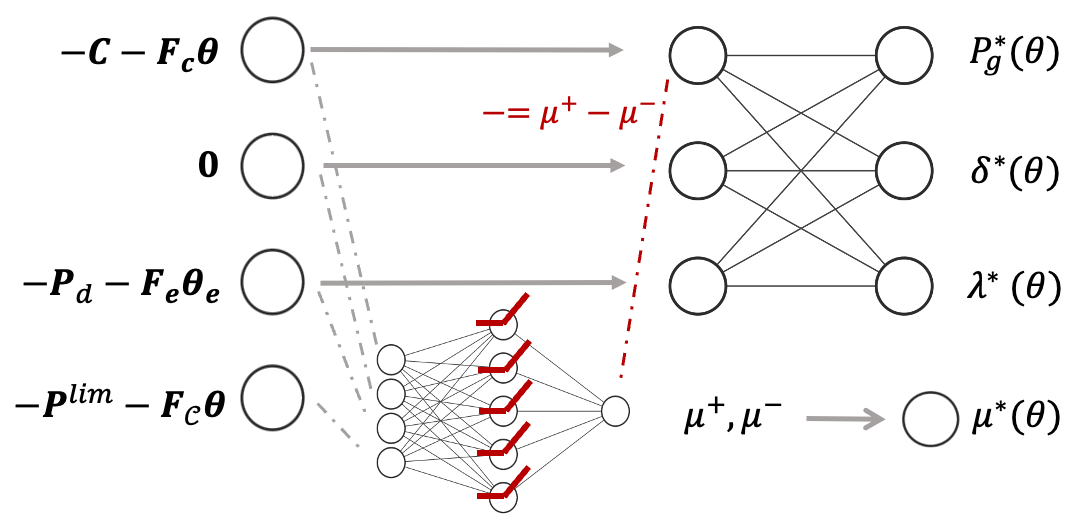}
		\caption{PSNN Architecture}
		\label{fig:architecture}
	\end{minipage}
\end{figure}

We use a partially supervised NN to predict the optimal solution through two sequential subnetworks: a shallow ReLU network defined in the form of (\ref{eq:deepNN}) of Section \ref{sec:NNs}, $\mu-$NN that calculates $\boldsymbol{\mu}^*$, and a separate NN that models $g(\bfm\mu^*;\bfm\theta)$ in (\ref{eq:linearsolver}). Figure \ref{fig:mupred} illustrates the subnetwork that predicts the optimal dual variables of the inequality constraints, $\boldsymbol{\mu}^*$. As discussed earlier, the weights of the first layer of this subnetwork model are fixed with $\nabla \boldsymbol{\mu}^*_{\mathcal{B}_i}$ vectors for all potential binding constraints $i \in \mathcal{C}$. Originally, each of these slope coefficient vectors have different dimensions, depending on which constraints are binding, but they are modified to have the same size by adding columns or rows containing zeros. The layer weights are obtained simply by stacking all the vectors vertically. The rest of the parameters, weights of the second layer, and biases, are determined via training. Finally, the second subnetwork is a mapping $g(\boldsymbol{\mu}^*;\bfm\theta) = [\textbf{P}_g^*,\boldsymbol{\delta}^*,\boldsymbol{\lambda}^*]$, which can be defined as a linear model (one layered NN) without training, by fixing its weights to $\textbf{W}=\textbf{J}^{-1}$ defined in eq \ref{eq:linearsolver}. Our PSNN model combines the two subnetworks in one NN flow as shown in Figure \ref{fig:architecture}. The $\mu-$NN takes $\textbf{C}, \textbf{P}_d,\textbf{P}^{+},\textbf{P}^{-}$ as inputs and produces $\boldsymbol{\mu}^*$. The second subnetwork takes the adjusted inputs to produce the optimal solutions.

\subsection{Test Results}

Our method was evaluated using IEEE-6, IEEE-30, IEEE-57 bus test systems implemented in MATPOWER package \cite{zimmerman2010matpower}. We used our PSNN model to predict optimal solutions to the DC-OPF problem, and compared them to solutions predicted by a NN model trained using classical methods (i.e., with paired input-output data representing particular values of system parameters and corresponding optimal solutions), and also to solutions generated analytically by the Gurobi solver (v10.0.3), which provides a baseline experimental control for comparing performance of the NN models. The PSNN models were trained using only 1000 data points, and validated and tested on datasets of 1000 and 4000 observations, respectively. To generate training and validation data, demand of all buses was set to 0.01, then one bus was chosen and its load incrementally increased from 0 to the maximum allowed by the system, and the procedure was repeated for one other bus. For the classically trained NN model, we generated training, validation and test sets of 5000, 1000 and 4000 observations, respectively, by setting the demand to the default load values in the Pandapower package \cite{thurner2018pandapower}, multiplying with a constant sampled from Uniform$(0.6,1.4)$ and solving the problem using Gurobi solver for each generated demand value \cite{zamzam2020learning,joswig2022opf}. All models were tested using two types of test datasets, generated to reflect realistic and extreme ranges of demand. Realistic demand data was generated by multiplying a base demand by a constant sampled from Uniform$(0.6,1.4)$. Extreme demand data was generated by fixing all load bus demands to 0.01, replacing one bus with a number sampled from $Uniform(0, \max)$ and repeating this step for all buses. All the experiments and training were done on a Mac Mini M2 (2023) on Python 3.9.18, using PyTorch 2.0.1.

The PSNN, classically-trained DNN, and Gurobi solver solutions were compared based on violation of the KKT optimality and feasibility conditions, namely (i) Stationarity eq. (KKT 1 $\textbf{P}_g, \bfm\delta$), (ii) Primal feasibility (KKT 2), (iii) Dual feasibility (KKT 3) and (iv) Complementary slackness (KKT 4), which are defined in Table \ref{tab:kkt} below. We calculated squared distance between the true optimal conditions and the model predictions as KKT 1, 2 and 4 are 0, and KKT 3 is non-positive at optimality.

\begin{table}
	\centering  
    \resizebox{\textwidth}{!}{
	\begin{tabular}{cccccc}
		KKT 1 $(\textbf{P}_g)$ & KKT 1 $(\boldsymbol{\delta})$ & KKT 2 & KKT 2 $(\leq)$ & KKT 3 & KKT 4 \\
		\midrule
		$\left(\frac{\partial L}{\partial \textbf{P}_g}\right)^2$ &
		$\left(\frac{\partial L}{\partial \boldsymbol{\delta}}\right)^2$ &
		$\left(\frac{\partial L}{\partial \boldsymbol{\lambda}}\right)^2$ &
		$\left[\max\left(0,\frac{\partial L}{\partial \boldsymbol{\mu}}\right)\right]^2$ &
		$\max(0,-\bfm\mu)$ &
		$\left(\boldsymbol{\mu}^*\frac{\partial L}{\partial \boldsymbol{\mu}}\right)^2$
	\end{tabular}}
	\caption{KKT Conditions and corresponding formulae}
	\label{tab:kkt}
\end{table}

%\begin{subequations}
%	\label{eq:kkt} 
%	\begin{align}
	%		\textbf{KKT-1} (\textbf{P}_g) &: \left(\frac{\partial L}{\partial \textbf{P}_g}\right)^2\\
	%		\textbf{KKT-1}  &: \left(\frac{\partial L}{\partial \boldsymbol{\theta}}\right)^2\\
	%		\textbf{KKT-2}  &: \left(\frac{\partial L}{\partial \boldsymbol{\lambda}}\right)^2\\
	%		\textbf{KKT-3}  &: \left[\max\left(0,\frac{\partial L}{\partial \boldsymbol{\mu}}\right)\right]^2\\
	%		\textbf{KKT-4}  &: \left(\boldsymbol{\mu}^*\frac{\partial L}{\partial \boldsymbol{\mu}}\right)^2
	%	\end{align}
%\end{subequations}

Table \ref{tab:KKTmean} reports mean squared KKT violations for the three models on the realistic test sets with local load characteristics. Note that this dataset has the same distribution characteristics as the dataset used to train the classical DNN model, but it is outside the training distribution for the PSNN model. The PSNN predictions satisfied all KKT conditions with less than 1E-10 MSE for 6- and 30-bus systems and less than 8E-8 MSE for 57-bus system, with relatively stable performance across the three test systems. On the other hand, KKT performance of the classically trained DNN models for 6- and 30-bus are close but is lower for 57-bus. KKT 1 MSE increased from 2.72E-06 and 7.04E-06 for 6-bus to 2.44E-04 and 5.9E-05 for 57-bus. Similar performance reductions can be seen for the KKT 2 and KKT 4 results, while the KKT 3 results were relatively stable with respect to the system size. Finally, Gurobi outperformed both NN based approaches in all cases. This is expected as the training data for PSNN and DNN were generated using the Gurobi solver, and the default settings of the package used to train the NN models allows for only 7 decimal points, so performance is limited by 1E-14 MSE for every prediction.

Table \ref{tab:KKTmeanExtreme} reports mean squared KKT violations for the three models on the test sets with extreme load characteristics. This test set reflects out-of-training distribution examples for both the classical DNN and PSNN models. Performance on the extreme datasets was close to the realistic demand results for PSNN, with 6- and 30-bus errors less than 1e-10 for all KKT measures. Performance reduced somewhat for the 57-bus system, where KKT 1 ($\bfm\delta$), KKT 3 and KKT 4 errors increased to 4.07E-08, 2.76E-08 and 2.62E-08 respectively. For this dataset, KKT violations of DNN performance is well above PSNN. Best performance was recorded for 30-bus system which produced at least $10^2$ times more squared error. Again, Gurobi outperformed the PSNN and DNN models in all cases as expected.

\begin{table}
	\centering
      \resizebox{\textwidth}{!}{
\begin{tabular}{lccccccc}
	\toprule
	& &KKT1-$\textbf{P}_g$ & KKT1-$\bfm\delta$ & KKT2 & KKT2 $(\leq)$ & KKT3 & KKT4\\
\midrule
	& \ 6bus& 4.79E-13 & 1.46E-10 & 1.44E-13 & 1.49E-14 & 3.47E-14 & 2.92E-14 \\
	% & 14bus & 9.02E-11 & 2.73E-08 & 1.28E-12 & 2.49E-08 & 2.68E-05 \\
	PSNN     & 30bus& 1.31E-13 & 1.17E-10 & 8.41E-14 & 0.00E+00 & 1.46E-11 & 1.11E-11 \\
	& 57bus& 6.68E-11 & 8.35E-08 & 2.54E-12 & 1.90E-10 & 1.37E-08 & 6.25E-08 \\
	
	\midrule
	& \ 6bus & 2.72E-06 & 7.04E-06 & 2.93E-05 & 9.70E-08 & 1.36E-08 & 1.01E-07 \\
	% & 14bus &  &  &  &  &  \\
	DNN & 30bus & 4.08E-07 & 1.86E-06 & 5.84E-05 & 0.00E+00 & 2.17E-09 & 6.00E-09 \\
	& 57bus & 2.44E-04 & 5.90E-05 & 5.21E-03 & 1.19E-08 & 3.27E-08 & 2.58E-06 \\
	\midrule
		&  \ 6bus & 1.52E-15 & 9.78E-28 & 2.64E-32 & 0.00E+00 & 0.00E+00 & 2.37E-16 \\
	Gurobi & 30bus & 6.80E-15 & 5.27E-12 & 1.23E-14 & 0.00E+00 & 0.00E+00 & 1.94E-19 \\
	&  57bus  & 4.01E-13 & 8.63E-10 & 5.96E-15 & 0.00E+00 & 0.00E+00 &  5.13E-19 \\
	\midrule
	\bottomrule
\end{tabular}}
			\caption{Mean Squared KKT Errors on the test sets with Local Perturbations}
	\label{tab:KKTmean}
\end{table}

\begin{table}[!ht]
	\centering

\resizebox{\textwidth}{!}{
\begin{tabular}{lccccccc}
	\toprule
	& &KKT1-$\textbf{P}_g$ & KKT1-$\bfm\delta$ & KKT2 & KKT2 $(\leq)$ & KKT3 & KKT4\\
	\midrule
	     & \ 6bus & 6.03E-13 & 1.61E-10 & 1.52E-13 & 1.40E-14 & 2.78E-14 & 2.64E-14 \\
 % & 14bus& 1.75E-09 & 1.24E-06 & 1.25E-10 & 1.94E-06 & 5.41E-03  \\
 PSNN   & 30bus & 1.64E-13 & 1.73E-10 & 2.04E-13 & 2.69E-12 & 2.38E-10 & 1.20E-10 \\
        & 57bus & 1.98E-11 & 4.07E-08 & 3.07E-12 & 2.76E-08 & 1.90E-09 & 2.62E-08 \\

      \midrule
    & \ 6bus & 7.15E-03 & 8.38E-01 & 3.26E+00 & 1.66E-02 & 5.89E-02 & 1.48E-03 \\
                   % & 14bus &  &  &  &  &  \\
   DNN & 30bus & 1.66E-02 & 1.54E-02 & 3.82E-02 & 4.76E-05 & 8.06E-07 & 2.32E-07 \\
       & 57bus & 7.41E+01 & 8.50E+01 & 3.28E+01 & 9.75E-04 & 4.03E-02 & 1.24E+00 \\
                                            \midrule
        & \ 6bus & 7.92E-14 & 3.00E-11 & 6.28E-16 & 0.00E+00 & 0.00E+00 & 2.37E-16 \\
        % & 14bus & 7.83E-13 & 4.00E-10 & 9.56E-15 & 0.00E+00 & 7.26E-16 \\
        Gurobi& 30bus & 6.80E-15 & 5.27E-12 & 1.23E-14 & 0.00E+00 & 0.00E+00 & 1.94E-19 \\
        & 57bus & 7.88E-13 & 1.36E-09 & 1.76E-14 & 0.00E+00 & 0.00E+00 & 1.30E-15 \\ 
  \bottomrule
	
\end{tabular}}
	\caption{Mean Squared KKT Errors on test sets generated for the Extreme Characteristics}
	\label{tab:KKTmeanExtreme}
\end{table}

\begin{table}\centering
  \resizebox{.9\textwidth}{!}{
	    \begin{tabular}{ccccccc}
		\toprule
		 & Case & Min  & 25\% & 50\% & 75\%  &  Max \\
		\midrule
		& 6bus & -7.24E-06 & -1.14E-06 & 7.42E-07 & 2.79E-06  & 9.38E-06\\
		Local & 30bus & -1.15E-06 & -3.02E-07 & -9.72E-08 & 1.02E-07 & 1.01E-06\\
		& 57bus & -1.15E-06 & -3.02E-07 & -9.72E-08 & 1.02E-07 & 1.01E-06\\
		\toprule
		\midrule
		& 6bus & -8.12E-06 & -7.27E-07 & 1.15E-06 & 2.81E-06 & 1.02E-05\\
		Extreme & 30bus & 2.56E-03 & 6.02E-03 & 7.12E-03 & 8.08E-03 & 9.80E-03\\
		& 57bus & -1.50E-03 & -5.06E-05 & -2.12E-05 & 2.74E-07 & 3.08E+01 \\
		\bottomrule
	\end{tabular}}
	\caption{$C(\textbf{P}^g) - \hat C(\textbf{P}^g)$ calculated for datasets constructed via local perturbations and reflecting extreme characteristics}
	\label{tab:OptimalityGap}
\end{table}

Table \ref{tab:OptimalityGap} reports the optimality gap between Gurobi solutions and our approaches. Here, we assume Gurobi as the reference point and calculated $C(\textbf{P}^g) - \hat C(\textbf{P}^g)$. For the test data, we report minimum, median and maximum optimality gaps along with 25th and 75th percentiles. The results show that the cost difference between our approaches and Gurobi is minimal for the vast majority of cases with median values around $\pm$1E-05. For the realistic dataset consisting of local demand perturbations, it can be seen that the median values calculated for PSNN are either very small or positive, indicating that our solutions cost close to or less than Gurobi solutions.  

% \begin{table}
	% 	\centering
	% 	\begin{tabular}{rrrr}
		% 		\toprule
		% 		6 Bus & 14 Bus & 30 Bus & 57 Bus \\ 
		% 		x63 & x235 & x2932 & x3543 \\ 
		% 		\bottomrule
		% 	\end{tabular}
	% 	\caption{Speed of NN predictions on the basis of Gurobi solver in finding solutions to 1000 examples}
	% 	\label{tab:speed}
	% \end{table}

\subsection{Uncertainty}
% The perks of using NNs for speed as intro. This allows us to use NNs for larger input parameters.

\begin{figure}
	\centering
	\includegraphics[scale=.3]{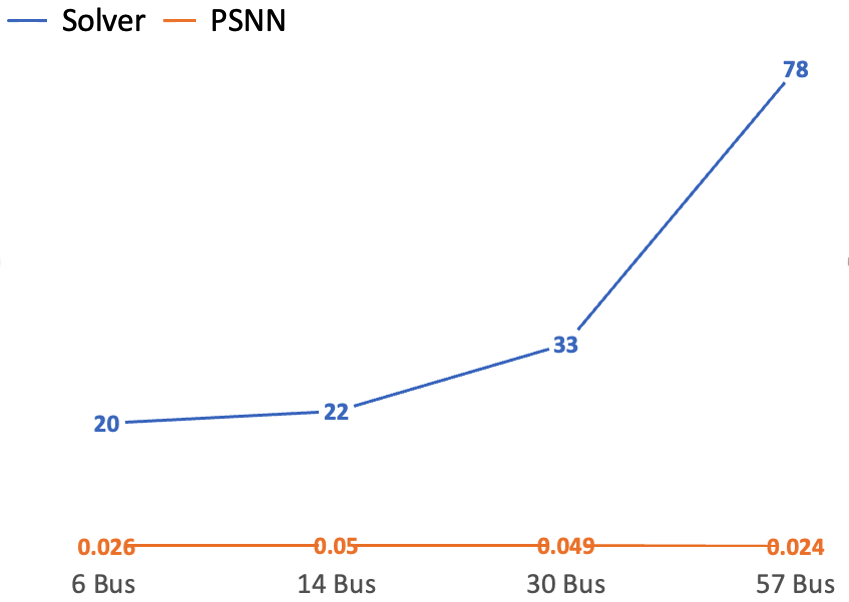}
	\caption{Calculation Speed Comparison of PSNN and Gurobi (Solver) over 1000 test problems}
	\label{fig:SSNN_speed}
\end{figure}

To compare the speed of our model against Gurobi, we randomly generated 1000 input parameters by applying local perturbations to the base load and solved the problem for all test systems. As shown in Figure \ref{fig:SSNN_speed} the solution time with Gurobi increased exponentially with system size. For the 6 bus system, it took 20 seconds to solve 1000 problems with Gurobi, which increased to 33 and 78 seconds for the 30- and 57-bus systems, respectively. On the other hand, solution time remained constant with the PSNN approach, ranging between 0.024 to 0.050 for any test system. 

The low computational cost and strong KKT performance of the PSNN model implies that it can be used to quickly generate large distributions of optimal and feasible solutions for the simulation and long term planning of energy systems with uncertain renewable resources, such as wind generation. To explore this we modified eq. \ref{eq:DC-OPF}b of the DC-OPF problem formulation as follows to allow variation:
\begin{equation}
	\textbf{B}\boldsymbol{\delta}-\textbf{P}_g-\textbf{P}_{ren} + \textbf{P}_d=0,
\end{equation}
where $\textbf{P}_{ren}\sim$ Exponential$(\lambda = 1.25)$. In the simulations, the problem was solved repeatedly for 500 random cases of $\textbf{P}_{ren}$ for each hour. A 1.5 unit generator limit was enforced by truncating the generated values above 1.5.

\begin{figure}
	\centering
	\begin{minipage}{.5\textwidth}
		\centering
		\includegraphics[width=\linewidth]{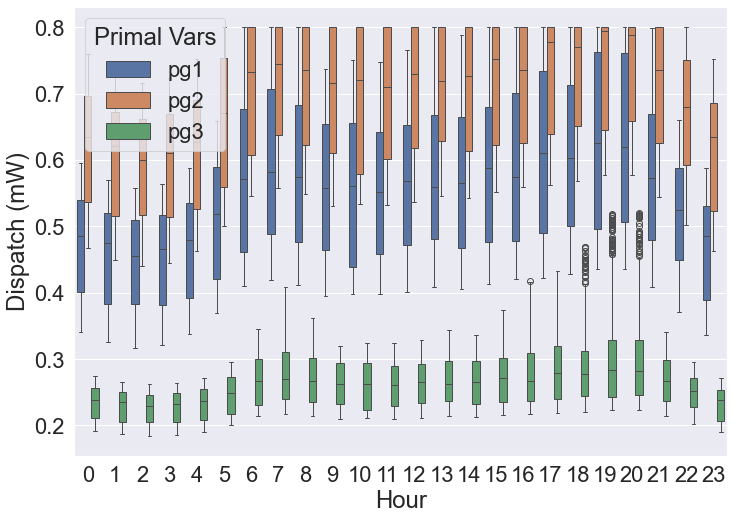}
		% \caption{Hourly distribution of optimal primal variables}
	\end{minipage}%
	\begin{minipage}{.5\textwidth}
		\centering
		\includegraphics[width=\linewidth]{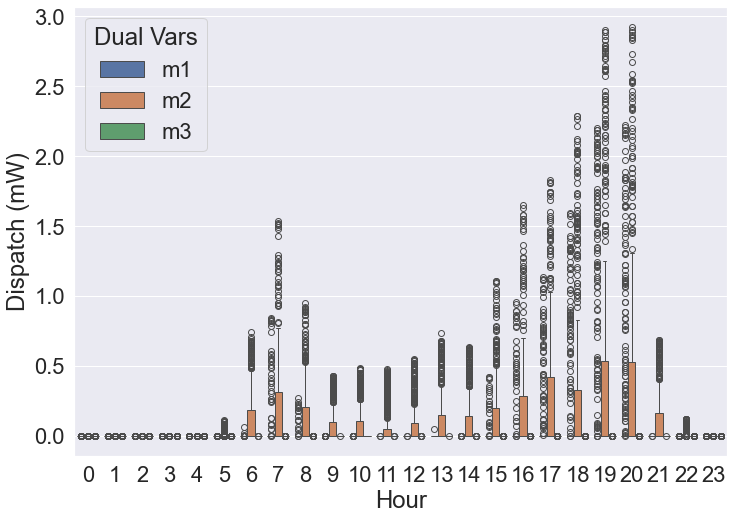}
	\end{minipage}
	\caption{Hourly distribution of optimal primal and dual prices}
	\label{fig:uncertain_boxplot}
\end{figure}

For this analysis,  default load values from the Pandapower package \cite{thurner2018pandapower} were taken as base demand. To reflect seasonality of hourly changing demand, we used historical demand data for Ontario, Canada \footnote{\url{https://www.ieso.ca/en/power-data/data-directory}} as a scaler. Specifically, an arbitrary day (May 11, 2024) was selected and the daily demand was divided by the maximum value throughout the day to act as scalers. Each bus value was multiplied by the constant for every hour to obtain 24 demand inputs. To calculate solutions for 1200 different inputs, the PSNN model ran for 0.03 seconds.

Figure \ref{fig:uncertain_boxplot} presents the distribution of the resulting 500 optimal generator dispatch for 24 hours using boxplots for the first three generators in the 30-bus system. All three generators have a capacity limit of 0.8 mW. As shown, the boxplots for the lowest cost Generator 2 are higher than other generators. During the day, Generator 2 often runs at full capacity and the dual variable of the associated capacity constraint ($P_{2} - 0.8\leq 0)$ is positive. During peak hours (5pm-9pm), Generator 1 also often reaches capacity with its dual variable increasing above zero, forcing the system to use Generator 3.

% Secondly, same test was performed for 13 days between April 29th to May 11th by running our model for 156000 inputs, which took 0.07 seconds to calculate. Figure \ref{fig:uncertain_lineplot} presents the calculated generator dispatches.

\section{Conclusion and Future Research} \label{sec:summary}
This paper proposes an explainable partially supervised NN  modeling approach that provides precise solutions for multiparametric QP optimization problems with linear constraints with generalizability power outside the training distribution. The model aligns the piecewise linear nature of NN architecture with the underlying mathematical structure of the optimization problem by deriving the model weights directly from the linear segments of the solution function in each critical region defined by the sets of binding constraints. To train the model, a cost efficient data generation approach was proposed that eliminates the need for using solvers in a loop to populate data. As a proof of concept, the PSNN modeling approach was applied to DC-OPF problems with upper and lower bound generator limits. PSNN models trained on only 1000 data points achieved higher KKT optimality and feasibility results than generic DNN models trained classically on five times as much data. The PSNN models generated optimal solutions to large input sets much faster than the Gurobi solver, with little sacrifice in KKT performance. As such, the models can rapidly create a distribution of optimal and feasible solutions to inputs sampled from empirical distributions of uncertain demand and renewable generator dispatches.

\section{Limitations and Future Work}\label{sec:ssnn_limitations}
\begin{figure}
	\centering
	\includegraphics[width=0.6\linewidth]{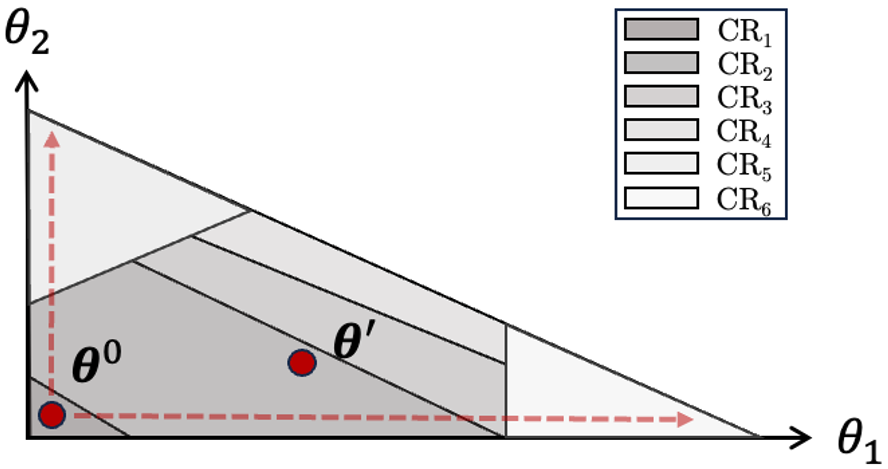}
	\caption{Illustration of critical regions of DC-OPF problem with line limits}
	\label{fig:lineLimits}
\end{figure}

\begin{figure}
	\centering
	\includegraphics[width=0.6\linewidth]{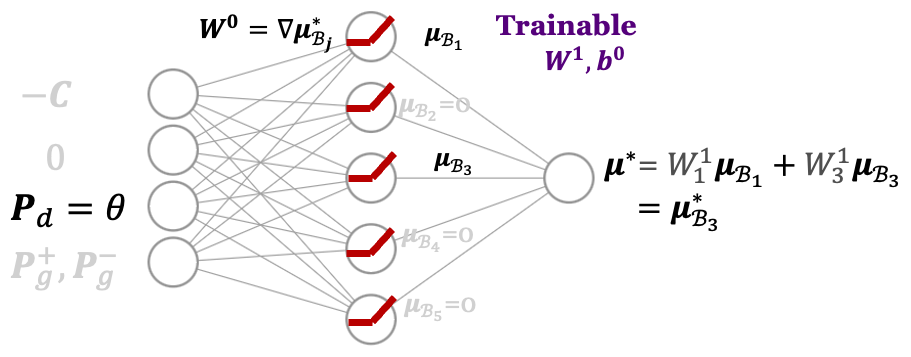}
	\caption{Illustration of Model Predictions}
	\label{fig:ssnn_linear_comb}
\end{figure}
This work has focused on the PSNN model architecture, its training, and providing proof-of-concept support. We acknowledge that further work is needed to develop an efficient discovery algorithm for identifying active constraint combinations that is applicable to a general class of  empirical applications. The PSNN methodology presented in this paper also has  limitations with respect to scalability. Our results indicated some declining performance with the larger DC-OPF systems. Further work is also needed to identify the difficulties in training PSNN models to scale up to larger systems.

The critical regions of the problem addressed in this paper, DC-OPF problem with box constraints, are relatively simple to discover. Specifically, fixing all parameters and moving along one axis is sufficient to discover all critical regions as the same pattern occurs in other dimensions. This is not a general pattern that occurs in other problems, such as DC-OPF with line limits, in which case moving along each axis would yield a different critical region as illustrated in Figure \ref{fig:lineLimits}.  Another limitation is the need to run a solver iteratively to discover the critical regions, by generating inputs by increasing the input parameter with a small increment, e.g., $\theta^{i} - \theta^{i-1}=0.01$. For larger systems and with a more comprehensive search pattern, computational cost of running solvers repeatedly can be very high.  However, this is not a limitation to our model itself but the search pattern that is employed as PSNN is capable of producing optimal solutions if the critical region has been discovered. Therefore, one can use PSNN for robust optimization instead where a smaller subset of all critical region is of interest. 

Second, while our models perform better than the alternative ones, for the 57-bus system, the training takes a large number of epochs to converge. The model was trained for 300,000 epochs until the training error reduced below 1E-10. This problem is due to finding the optimal trainable model weights and biases that can choose the correct slope from the first layer of the NN. To illustrate this point mathematically, consider the model
\begin{align}
	f_\mu(-\textbf{B}-\bfm\theta) = \textbf{W}^1\sigma(\textbf{W}^0\textbf{x}+b^0).
\end{align}
As the slopes of $\mu(\bfm\theta)$ for all critical regions are contained in the first layer weights, the output of the first layer is all candidate slopes. After the ReLU activation function is applied, the negative predictions are replaced with 0, which yields
\begin{equation}
	\textbf{h}^0 = \sigma(\textbf{W}^0(-\textbf{B}-\bfm\theta)) = \sigma( [\bfm\mu^+_{\mathcal{B}_1},\bfm\mu^+_{\mathcal{B}_2},\dots,\bfm\mu^+_{\mathcal{B}_k}]^T),
\end{equation}
for $k$ number of critical regions, where $\bfm\mu^+_{\mathcal{B}_j} = \bfm\mu^*_{\mathcal{B}_j}$ if the result is nonzero and 0 otherwise. Then, the rest are linearly combined with the second layer weights, as illustrated in Figure \ref{fig:ssnn_linear_comb}, where the optimal solution, $\bfm\mu^*=\bfm\mu^*_{\mathcal{B}_3}$, is obtained by linearly combining $\bfm\mu^*_{\mathcal{B}_1}$ and $\bfm\mu^*_{\mathcal{B}_3}$ when $\bfm\theta\in CR_3$. Finding such model parameters that would apply all critical regions is a challenging task that will be addressed in future work.

% While our approach shows a more direct way to integrate DNNs to mp-QP solutions, it suffers from challenges in training, which makes it harder to reduce the training error to zero even though it is possible to do so. As a result, it is harder to scale the model for higher dimensional problems, which is also a problem of classical approaches to NN. 

% Related to this point, during the discovery steps, the parameters are fixed to a small number and one parameter, e.g. $\theta_j$, is chosen randomly to be increased from 0 to the largest feasible number. This process discovers lower and upper bounds of critical regions for slopes $\nabla \bfm\mu_{\mathcal{B}_i}^*$ such that when $\theta_i \in [CR^-_j,CR^+_j]$, $\nabla \bfm\mu_{\mathcal{B}_j}^*$ produces the optimal solutions. To populate data, we repeated data sampling along another axis assuming the CRs are bounded by the same intervals, $\theta_{i'}\in [CR^-_j,CR^+_j]$ and removed observations that violate KKT optimality conditions. 

\bibliography{refs.bib}
\bibliographystyle{unsrt}

\appendix

\section{Supplemental material}

\subsection{Neural Networks as Piecewise Linear Functions}

\label{sec:NN_PWL}
The proposed PSNN approach relies on the observation that shallow NNs with ReLU activation function fall under the class of PWL functions whose slope and intersection parameters are estimated from data. As a PWL function, a shallow NN with $d_h$ neurons can approximate any continuous function by fitting $d_h+1$ linear segments. Thus, as $d_h$ increases, the estimation of the function improves, given that the model is trained with a sufficiently large dataset. To support later discussion, this section formulates NN as PWL functions.

% However, more complex models, i.e., $d_h>n$, cannot guarantee to estimate a piecewise linear target function with $n$ connected lines as it can zero the training error in an infinite variety of ways. 

A shallow NN with $d_h$ neurons and no activation function in the output layer can be written as
\begin{equation}
	\label{eq:BG_shallowNN}
	f(\textbf{x};\varphi) = \textbf{W}^{1T}\sigma(\textbf{W}^{0T}\textbf{x}^T+\textbf{b}^0) + \textbf{b}^1,
\end{equation}
where $\varphi = [\textbf{W}^0,\textbf{W}^1,\textbf{b}^0,\textbf{b}^1]$ is the parameter set, $\textbf{W}^0 \in \mathbb{R}^{d_h\times d_x}$, $\textbf{W}^1 \in \mathbb{R}^{d_y\times d_h}$, $\textbf{b}^0 \in \mathbb{R}^{d_h}$, $\textbf{b}^1 \in \mathbb{R}^{d_y}$, $\sigma(\textbf{z})=ReLU(\textbf{z})$, $\textbf{x} \in \mathbb{R}^{d_d\times d_x}$, $d_x, d_y$ are input and output dimensions of the NN, respectively, and $d_d$ is the number of observations. 

For the case when $d_y=1$, for a given observation, $\textbf{x}_i \in \mathbb{R}^{d_x}$, hidden layer weights $\textbf{W}^0_j \in \mathbb{R}^{d_x}, \textbf{W}^1_j \in \mathbb{R}$ and scalar $b_j$ for $j\in [1,\dots,d_h]$, the model can be rewritten as
\begin{align}
	f(\textbf{x}_i;\varphi) = \sum_{j=1}^{d_h} \textbf{W}^1_{j}(\textbf{W}^0_{j}\textbf{x}_i+b^0_j)^+  + b^1,
\end{align}
where $(z)^+= z$ if $z>0$ and 0 otherwise. 

It can be seen that the above function is a weighted summation of hidden layer output, which generates a different linear function within a region of $\textbf{x}_i$ depending on whether one or more nodes are active, i.e., $(\textbf{W}^0_{j}\textbf{x}_i+b^0_j)>0$ in that region. Expressed mathematically, let $\mathcal{B}_k$ be the set of indices of nodes that are active for a given $\textbf{x}$, i.e., $\mathcal{B}_k =\{j \ | \textbf{W}^0_{j}\textbf{x}+b^0_j >0 \}$. Also let $\mathcal{R}_k \subseteq \mathcal{D}$ be a subdomain of all $x$ in which the nodes with indices in $\mathcal{B}_k$ are active, i.e., $\mathcal{R}_k =\{\textbf{x} \ | \ \textbf{W}^0_{j}\textbf{x}+b^0_j >0, j\in \mathcal{B}_k\}$. 

Notice that for a compact domain of $\textbf{x} \in \mathcal{D}$, there is a finite number $K$ such that $\bigcup_{k=1}^K \mathcal{R}_k = \mathcal{D}$ and $\mathcal{R}_k \cap \mathcal{R}_{k'} = \emptyset$ for $k\neq k'$. Therefore we can write,

\begin{align}
	f&(\textbf{x};\theta) = b^1 + \notag \\
	& \begin{cases}
		0 & \text{ for }  \textbf{x} \in \mathcal{R}_0\\
		\sum_{j \in \mathcal{B}_1} \textbf{W}^1_{j}\textbf{W}^0_{j}x_i+ \textbf{W}^1_{j}b^0_j  & \text{ for } \textbf{x} \in \mathcal{R}_1 \\\
		\vdots & \vdots\\
		\sum_{j \in \mathcal{B}_{K-1}} \textbf{W}^1_{j}\textbf{W}^0_{j}x_i+ \textbf{W}^1_{j}b^0_j   & \text{ for } \textbf{x} \in \mathcal{R}_{K-1}\\
		\sum_{j \in \mathcal{B}_K} \textbf{W}^1_{j}\textbf{W}^0_{j}x_i+ \textbf{W}^1_{j}b^0_j   & \text{ for } \textbf{x} \in \mathcal{R}_K\\
	\end{cases}.
\end{align}
For each $\mathcal{R}_k$, the NN can generate a different linear function. Also, the maximum number of linear functions that a model can generate is limited by the hidden layer size, $K=d_h+1$. 

By  substituting $\textbf{M}_i = \sum_{j \in \mathcal{B}_i} \textbf{W}^1_{j}\textbf{W}^0_{j}\textbf{x}$, and $ \textbf{N}_i = \sum_{j \in \mathcal{B}_i}\textbf{W}^1_{j}b^0_j $, one can obtain the general form of a PWL function,
\begin{align}
	f(x_i,\theta) = \begin{cases}
		\textbf{N}_0 & \text{ for }  \textbf{x} \in  \mathcal{R}_0\\
		\textbf{M}_1\textbf{x}+ \textbf{N}_1  & \text{ for } \textbf{x} \in  \mathcal{R}_1 \\\
		\vdots & \vdots\\
		\textbf{M}_{K-1} \textbf{x}+ \textbf{N}_{K-1}  & \text{ for } \textbf{x} \in  \mathcal{R}_{K-1}\\
		\textbf{M}_K\textbf{x}+ \textbf{N}_K  & \text{ for } \textbf{x} \in  \mathcal{R}_K\\
	\end{cases}.
\end{align}
% As can be seen above, the number of cases in the function above is limited by the number of hidden neurons, $d_h$, that can be activated.

\begin{figure}
	\centering
	\includegraphics[scale=.4]{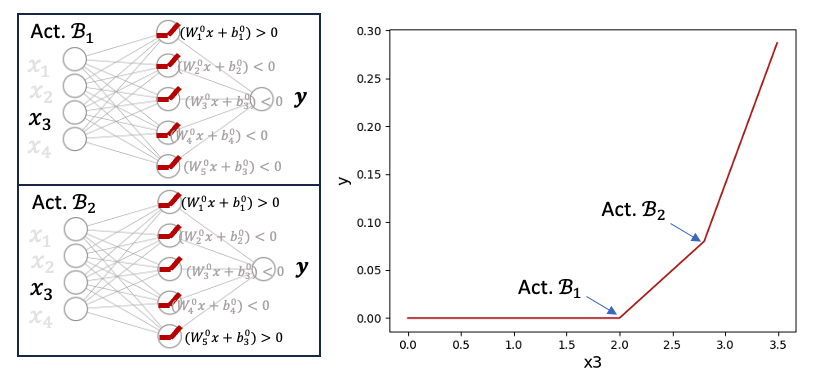}
	\caption{Illustration of NN as a Piecewise Linear Function}
	\label{fig:plf}
\end{figure}
Figure \ref{fig:plf} illustrates how each activated neuron changes the slope of the piecewise linear function represented by a NN. The points donated as $Act. \mathcal{B}_1$ and $Act. \mathcal{B}_2$ show where different sets of binding constraints are activated and NN diagrams on the left show activated neurons at respective points.  As illustrated, for $x_3<2$, none of the hidden neurons have a positive value as the $\textbf{W}_j^0x+b_j<0$ for all $j$. In the interval $2 \leq x_3 <4$, the set of active neurons is $\mathcal{B}_1 = \{1\}$ as only $\textbf{W}_1^0x+b_1>0$ and the output increases with a certain slope as $x_3$ increases. After $x_3\geq 4$, another slope is also activated and the set of active slopes are $\mathcal{B}_2 = \{1,5\}$, which updates the change of the output value with respect to $x_3$.

If the target function is not PWL, but an arbitrary nonlinear continuous function defined in a bounded interval, $x_i\in [a,b]$, then the NN can be used to approximate the function, and the prediction accuracy increases as $d_h \rightarrow \infty$, in accordance with the Universal Approximation Theorem (UAT) \cite{hornik1989multilayer}.

\subsection{Derivation of Equality Constrained Problem Solver}
\label{app:linear_derivation}
The gradient of the Lagrangian function (\ref{eq:lagr_eq}) with respect to $\textbf{x},\boldsymbol{\lambda}$ must be zero at optimality as follows:
\begin{subequations}
	\begin{align}
		\frac{\partial L}{\partial \textbf{x}} &= 2\textbf{Q}\textbf{x}+\textbf{C} + \bfm\theta_c - \boldsymbol{\lambda}^T\textbf{A}_e=0\\
		\frac{\partial L}{\partial \boldsymbol \lambda} &= \textbf{b}_e+\bfm\theta_c-\textbf{A}_e\textbf{x} =0.
	\end{align}
	\label{eq:derivatives}
\end{subequations}
Since, the parameter $ \bfm\theta $ in (\ref{eq:derivatives}b) can change due to the variability of the studied problem, such as electricity demand in power system application and non-dispatchable supplies such as wind and solar generators, our goal is to train a model that predicts the optimal solution of (\ref{eq:ineq}a)-(\ref{eq:ineq}b) corresponding to an arbitrary $\bfm\theta$,  while the remaining system parameters $\textbf{Q}, \textbf{C},\textbf{A}_e,\textbf{F}_c,\textbf{A}_e,\textbf{A}_\mathcal{C}$ are considered constant. Therefore, equations (\ref{eq:derivatives}a) and (\ref{eq:derivatives}b) can be written in matrix form as,
% rearranged into a one-to-one and multivariate function, which can be written as,
\begin{align}
	\begin{bmatrix}
		2\textbf{Q} & -\textbf{A}_e^T\\
		-\textbf{A}_e & 0
	\end{bmatrix}
	\begin{bmatrix}
		\textbf{x}\\ \boldsymbol{\lambda}
	\end{bmatrix} 
	=
	\begin{bmatrix}
		-\textbf{C}-\bfm\theta_c \\ -\textbf{b}_e-\bfm\theta_c
	\end{bmatrix},\label{eq:linear}
\end{align}
where the optimal solution to (\ref{eq:ineq}a)-(\ref{eq:ineq}b) for a given $\bfm\theta$ can be obtained by solving (\ref{eq:linear}). To this end, the $ \textbf{J} $ matrix can be defined as below, whose inverse is used to find the optimal solution as
\begin{align}
	\textbf{J}^{-1}
	\begin{bmatrix}
		-\textbf{C}-\bfm\theta_c \\ -\textbf{b}_e-\bfm\theta_c
	\end{bmatrix}
	=
	\begin{bmatrix}
		\textbf{x}^*\\ \boldsymbol{\lambda}^*
	\end{bmatrix}, \quad \text{for} \quad \textbf{J} = \begin{bmatrix}
		2\textbf{Q} & -\textbf{A}_e^T\\
		-\textbf{A}_e & 0
	\end{bmatrix}.
	\label{eq:linearinverse}
\end{align}
%			\subsection{Equality Constrained Solution Function}
%	  The mapping in (\ref{eq:linear}) can be expressed as a linear function, $g^{-1}:[\textbf{x},\boldsymbol{\lambda}]\rightarrow [\textbf{e}_1-\textbf{C},\textbf{e}_2-\textbf{b}_1]$, whose inverse $g:[\textbf{C},\textbf{b}_1]\rightarrow [\textbf{x}^*,\boldsymbol{\lambda}^*]$, is also a function, which can be written as,		
From the above equation \ref{eq:linearinverse},  the function $ g(\textbf{C},\textbf{b}_e) $ can be obtained as
\begin{align}
	g(\textbf{C},\textbf{b}_e+\bfm\theta_c) = \begin{bmatrix}
		2\textbf{Q} & -\textbf{A}_e^T\\
		-\textbf{A}_e & 0
	\end{bmatrix}^{-1}
	\begin{bmatrix}
		-\textbf{C}-\bfm\theta_c \\ -\textbf{b}_e-\bfm\theta_c
	\end{bmatrix}.
\end{align}

\subsection{Derivation of Lagrangian Derivatives for DC-OPF}
\label{app:dc_derivation}
\begin{align}
	L(\textbf{P}_g,\boldsymbol{\theta},\boldsymbol{\lambda},\boldsymbol{\mu})&=
	\textbf{P}_g^T\textbf{Q}\textbf{P}_g+(\textbf{C}^T+\bfm\theta_c)\textbf{P}_g+\textbf{C}_0 + \notag \\ &\quad\boldsymbol{\lambda}(\textbf{P}_d-\textbf{P}_g-\textbf{B}\boldsymbol{\delta})+ \boldsymbol{\mu}^+(\textbf{P}_g-\textbf{P}_g^{+}) + \boldsymbol{\mu}^-(\textbf{P}_g^{-}-\textbf{P}_g).
\end{align}
% where $\boldsymbol{\mu}^T = [\boldsymbol{\mu}^+,\boldsymbol{\mu}^-]$. 
Following the same steps in Section \ref{sec:method}, our model is derived using the derivatives of the Lagrangian function  of the DC-OPF problem can be derived as below

The partial derivatives of the above function is
\begin{subequations}
	\begin{align}
		\frac{\partial L}{\partial \textbf{P}_g} &= 2\textbf{Q}\textbf{P}_g+\textbf{C}_1+\boldsymbol{\lambda}+\boldsymbol{\mu}=0,\\
		\frac{\partial L}{\partial \boldsymbol{\theta}} &= \textbf{B}^T\lambda = 0,\\
		\frac{\partial L}{\partial \boldsymbol{\lambda}} &= \textbf{P}_d-\textbf{P}_g-\textbf{B}\boldsymbol{\delta}=0,\\
		\frac{\partial L}{\partial \boldsymbol{\mu}} &= \textbf{P}_g - \textbf{P}_g^{max} = 0.
		% \frac{\partial L}{\partial \boldsymbol{\mu}^-} &= \textbf{P}_g^{min} - \textbf{P}_g = 0.
	\end{align}
\end{subequations}

\subsection{Further Details on Training of PSNN}
\label{app:training}
All models were compared on two types of datasets, realistic and extreme. The former was constructed by multiplying a base load with a factor drawn from Uniform(0.6,1.4). The latter was constructed by choosing one $\theta_i$ and sampling $m$ values from Uniform(0,\textit{max}) while fixing the remaining $\theta_j$ to 0.01 and repeating this for all parameters. For this method, \textit{max} was set to the sum of all generator capacity limits. To generate $M$ number of data points, $m$ is set to $\lceil M/n_l\rceil$ where $n_l$ is the number of dimensions. If more than $M$ points were sampled, the excessive amount was removed.

While PSNN was tested on these datasets, the training set is constructed using the critical regions $[\textbf{CR}_i^-,\textbf{CR}_i^+]$ detected via the discovery step. To discover sets of binding constraints, all parameters were set to 0.01 and an arbitrary $\theta_i'$ was incrementally increased from 0 to the maximum number the system allows. The discovery was carried out on one dimension only, so this approach would not detect if other sets of binding constraints could be found by repeating this step with other $\theta_i$. 

To construct the training set for PSNN, different values of $\theta_i'$ were sampled from each $[\textbf{CR}_i^-,\textbf{CR}_i^+]$ and solved by expanding the coefficient matrix, $\textbf{J}_\mathcal{B_i}$ with the corresponding binding constraints, running the calculations in eq. \ref{eq:linearinverse_ineq} and checking if the solutions satisfy the KKT conditions. This step was repeated using one other $\theta_j$ and fixing all other parameters, assuming the same $[\textbf{CR}_i^-,\textbf{CR}_i^+]$ would apply to this range. If the calculated solution does not satisfy KKT, another number was sampled from the same interval. 

The problem of critical regions that do not generalize to other dimensions could be solved by calculating optimal solutions using all discovered $\nabla \bfm\mu^*$ and choosing the ones that satisfies KKT conditions. For this study, however, our models were trained using the above mentioned approach.

\end{document}